\documentclass[twoside,final,11pt]{entics}
\input{localnom.sty}
\usepackage{enticsmacro}
\usepackage{graphicx}
\usepackage[all]{xy}

\def\conf{MFPS 2026}
\volume{6}
\def\volu{6}
\def\papno{1}
\def\mydoi{19496}
\def\lastname{Bayeh, Fu, Selinger}
\def\runtitle{The category of nominal sets is locally monoidal closed}
\def\copynames{F. Bayeh, P. Fu, P. Selinger}

\begin{document}
\begin{frontmatter}
  \title{The Category of Nominal Sets is Locally Monoidal Closed}
  \author{Fahimeh Bayeh\thanksref{a}\thanksref{fahimeh.bayeh@dal.ca}}
  \author{Peng Fu\thanksref{b}\thanksref{pfu@cse.sc.edu}}
  \author{Peter Selinger\thanksref{a}\thanksref{selinger@dal.ca}}
  \address[a]{Department of Mathematics and Statistics\\Dalhousie
    University\\Halifax, Nova Scotia, Canada}
  \address[b]{Department of Computer Science and Engineering\\University
    of South Carolina\\Columbia, South Carolina, U.S.A.}
  \thanks[fahimeh.bayeh@dal.ca]{Email: \href{mailto:fahimeh.bayeh@dal.ca} 
    {\texttt{\normalshape fahimeh.bayeh@dal.ca}}}
  \thanks[pfu@cse.sc.edu]{Email: \href{mailto:pfu@cse.sc.edu} 
    {\texttt{\normalshape pfu@cse.sc.edu}}}
  \thanks[selinger@dal.ca]{Email: \href{mailto:selinger@dal.ca} 
    {\texttt{\normalshape selinger@dal.ca}}}
  
  \begin{abstract}
    The category $\Nom$ of nominal sets was proposed by Pitts and Gabbay
    as a setting for the semantics of abstract syntax with variable
    bindings. It is well-known that $\Nom$ is a topos, also known as the
    Schanuel topos, and in particular it follows that $\Nom$ is
    \emph{locally closed}, i.e., every slice category $\Nom/X$ is
    cartesian-closed. In this paper, we show that $\Nom$ has a much
    stronger property: every monoidal (closed) structure on $\Nom$
    induces a monoidal (closed) structure on all slice categories
    $\Nom/X$. One monoidal closed structure of particular interest on
    $\Nom$ is the \emph{separated product} $A\ast B$, whose right adjoint
    $A\magicwand B$ is called the \emph{separated function space}. In
    particular, it follows that $\Nom$ is locally separatedly closed.
  \end{abstract}
  
\end{frontmatter}

\section{Introduction}

The category $\Nom$ of nominal sets was proposed by Pitts and Gabbay
{\cite{GabbayPitts1999,GabbayPitts2002}} as a setting for the semantics of abstract
syntax with variable bindings. This category has a wealth of
interesting structure. Among other things, it has all limits and
colimits, it is cartesian closed, and in fact a topos (in this
context, it is also known as the \emph{Schanuel topos}). But the
feature that makes $\Nom$ uniquely interesting is the presence of
\emph{name abstraction}: for every object $X$ (whose elements we may
think of as the \emph{terms} of some syntax), there is an object
$[\aA]X$, whose elements represent \emph{abstractions}, which we can
think of as terms $a.t$ with a bound variable $a$. Crucially, the
elements of $[\aA]X$ are automatically defined up to
$\alpha$-equivalence, i.e., up to renaming of bound variables. This
makes $\Nom$ suitable as a category in which to interpret abstract
data types with binders. For example, the nominal set $T$ of lambda terms can be defined as the initial fixed point of
the equation $T = \aA + T\times T + [\aA]T$, reflecting the fact that
the set of terms is freely generated so that every term is either a
variable ($\aA$), an application ($T\times T$), or a lambda
abstraction ($[\aA]T$). The fact that this simple definition
automatically takes care of bound variables and $\alpha$-equivalence
is what makes nominal sets so useful in this context.

Theorem provers such as Rocq, Agda, and Lean are based on dependent
type theory {\cite{rocq,agda,lean}}. Seely showed that an appropriate
categorical structure for giving semantics of dependent type theories
is a locally cartesian closed category, i.e., a category $\Cc$ such
that all of its slice categories $\Cc/X$ are cartesian closed
{\cite{Seely1984}} (see also \cite{Clairambault2014}).  Since $\Nom$
is locally cartesian closed, one might imagine that it can support the
semantics of a dependent type theory in which one can also reason
about data types with binders \cite{Pitts2013}. However, defining such
a type theory has so far proved to be elusive; some steps in this
direction were taken by
{\cite{Cheney2012,Pitts2015,Fairweather_et_al2015,VanMuylder_et_al2025}}. One
of several issues that need to be addressed is that in such a
dependent type theory, it is not sufficient to have a global object of
name abstractions; instead, we must consider kinding judgements such
as $\Gamma\vdash [\aA]X : \type$, where $X$ depends on variables
defined in the context $\Gamma$.  As is usual in dependent type
theory, we will interpret a kinding judgement $\Gamma \vdash X :
\type$ as a map $X\to\Gamma$, i.e., an object in $\Nom/\Gamma$. Given
such an object, we then need to define the interpretation of $[\aA]X$,
also as an object in $\Nom/\Gamma$.  Therefore, not just $\Nom$
itself, but also all of its slice categories $\Nom/\Gamma$ should
support name abstractions and other related constructs.

In this paper, we begin to address this issue by showing that the
required local structure indeed exists in $\Nom$. In fact, we prove
something more general. Sch\"opp {\cite{Schopp2006}} showed that name
abstraction $[\aA]X$ arises as a special case of a more general
monoidal closed structure on $\Nom$, given by the so-called
\emph{separated product} $\ast$ and its right adjoint $\magicwand$,
called the \emph{separated function space}. The definition of
separated function spaces was later simplified by Clouston
{\cite{Clouston2013}}. We show that $\Nom$ is \emph{locally
  separatedly closed}, i.e., there is a version of $\ast$ and
$\magicwand$ on every slice category $\Nom/\Gamma$.  However, in
proving this, we found something that surprised us: one can
show much more generally that \emph{every} monoidal (respectively,
monoidal closed) structure on $\Nom$ lifts to a monoidal
(respectively, monoidal closed) structure on $\Nom/\Gamma$.

This is surprising because it seems to be very rare. In generic
categories, not only is it not usually possible to transport monoidal
structure from the base category to its slice categories, but there is
not even any canonical notion of what it would mean to do so. In other
words, there is no obvious way in which one can say that a monoidal
structure on a slice category is ``the same as'' that on the base
category --- and consequently, there is no accepted definition of a
``locally monoidal'' (or ``locally monoidal closed'') category.
It seems to be a particular feature of $\Nom$ that there is such a
notion, and moreover all such structures are local.

The structure of the paper is as follows: Section~2 covers the basics
of nominal sets. In Section~3, we discuss how one can extend the
$\Pi$-action from finitely supported permutations to all
permutations. In Section~4, we define how every bijection between atom
sets induces a functor between the category of nominal sets over those
atom sets.  In Section~5, we give the definition of nominal families
over a nominal set $X$ and prove that the category of nominal families
over $X$ and the slice category over $X$ are equivalent. Section~6
covers the monoidal structure of the slice category $\Nom/X$, while
Section~7 covers the monoidal closed structure of $\Nom/X$.

\section{Background}

In this section, we will review some of the basic notions of nominal
set theory. For a much more detailed introduction, see \cite{Pitts2013}.

\subsection{Names and permutations}

We fix a countably infinite set $\aA$, whose elements will be called
\emph{atoms} or \emph{names}. Atoms serve as unique identifiers. Their
two main properties are that we can always find a previously unused
(``fresh'') atom, and that they can be compared for equality.

\begin{xdefinition}
  A \emph{finitely supported permutation} of $\aA$ is a permutation
  $\pi:\aA\to\aA$ such that the set $\set{a\in \aA\mid \: \pi(a) \neq
    a }$ is finite. Let $\PIA$ be the group of these finitely
  supported permutations. We drop the subscript when it is clear from
  the context. The group operation is the composition of permutations,
  which we denote by $\pi_1\circ\pi_2$ or simply $\pi_1\pi_2$. The
  identity permutation is denoted by $\id$.
\end{xdefinition}

\begin{xdefinition}
A $\PI$-set is a pair $(X,\bul)$, where $X$ is a set and $\bul$ is an action of $\PI$ on $X$. We often just write $X$ for a $\PI$-set when the action is clear from the context.
\end{xdefinition}

\begin{xdefinition}
Let $X$ and $Y$ be $\PI$-sets. A function $f:X\to Y$ is \emph{equivariant} if for all  $\pi\in \Pi$ and for all $x\in X, \pi\bul \left( f(x) \right)=f\left( \pi\bul x \right)$.
\end{xdefinition}

\begin{xdefinition}
  Let $X$ be a $\PI$-set. Then the action of $\PI$ on $\powerset(X)$,
  the powerset of $X$, is defined by $\pi \bul S := \set{ \pi a
    \mid a\in S }$. This means that $\powerset(X)$ is a $\PI$-set.
\end{xdefinition}

\begin{xdefinition}\label{def:function-pi-action}
Let $X$ and $Y$ be $\PI$-sets and let $\pi \in \PI$. The action of
$\pi$ on a function $f:X\to Y$ is defined by $(\pi \bul f)(x) = \pi
\bul \left(f \left( \pi\inv \bul x\right) \right)$.
\end{xdefinition}

\subsection{Support}

We write $X\setminus Y$ to denote the set difference, so $X\setminus
Y=\set{x\in X\mid x\not\in Y}$.

\begin{xdefinition}
  Let $X$ be a $\PI$-set and let $x\in X$. We say that the set
  $S\subseteq \aA$ is a \emph{support} for $x$ if for all $a,b\in
  \aA\setminus S$, $(a\,b) \bul x=x$. Given a $\PI$-set $X$,
  $x\in X$ is \emph{finitely supported} if there exists some finite
  support $S$ for it. One can show that for any finitely supported
  $x$, there exists a smallest finite support, and we denote it by
  $\supp(x)\in P(\aA)$ \cite{Pitts2013}.
\end{xdefinition}

\begin{xdefinition}
  A $\PI$-set $X$ is called a \emph{nominal set} if all $x\in X$ are
  finitely supported.
\end{xdefinition}

\noindent Every set is nominal with the trivial action (``discrete''
nominal set). One non-trivial example of a nominal set is $\aA$
itself, with the natural action $\pi\bul a=\pi(a)$. Moreover, if $X,Y$
are nominal sets, then the set of finitely supported functions
$f : X \to Y$ is also a nominal set, with the action as in
Definition~\ref{def:function-pi-action}. Also, if $f:X\to Y$ is
finitely supported, then for all $x\in X$, we have
$\supp(f(x)) \subseteq \supp(f) \cup \supp(x)$. The function $f$ is
equivariant if and only if $\supp(f)$ is empty; in this case, we have
$\supp(f(x)) \subseteq \supp(x)$. We will discuss other examples of
nominal sets below.

\subsection{Freshness}

Let $X$ and $Y$ be nominal sets and consider $x\in X$ and $y\in Y$. We
say $x$ is \emph{fresh for} $y$, denoted $x\# y$, if they have
disjoint supports, i.e., if $\supp(x) \cap \supp(y) = \emptyset$.
Note that for all $\pi \in \PI$, we have $x\# y \Rightarrow (\pi
\bul x) \# (\pi \bul y)$. In other words, the freshness relation
is equivariant.

Now consider a finitely supported function $P:\aA\to \Bool$, which we
think of as a property of elements of $\aA$. One can show that $P(a)$
holds for some $a\# P$ if and only if $P(a)$ holds for all $a\# P$ if
and only if $P(a)$ holds for all but finitely many $a\in\aA$. This
motivates the \emph{freshness quantifier}: we are justified in saying
``for fresh $a$'' instead of ``for all fresh $a$'' or ``for some
fresh $a$''.

\subsection{The category $\Nom$}

Fix an atom set $\aA$. The category $\NomA$ has nominal sets as objects
and equivariant functions as morphisms. We omit the subscript when
it is clear from the context.

The category $\Nom$ has arbitrary coproducts, which are given by
disjoint union with the obvious $\PI$-action. It also has finite
products, where the $\PI$-action on $X\times Y$ is given by
$\pi\bul(x,y) = (\pi\bul x, \pi\bul y)$. It also has infinite
products, but they are not the same as in $\Set$: we must consider the
subset of the infinite cartesian product that consists of finitely
supported tuples. The category $\Nom$ is also cartesian closed, where
the nominal function space $\NFS{X}{Y}$ is the set of all finitely
supported functions $f:X\to Y$. In addition, $\Nom$ has many other
interesting properties; for example, it is a topos.

\subsection{Separated closure}

In this paper, we are particularly interested in monoidal closed
structures on $\Nom$. Besides the cartesian closed structure, there is
one other monoidal closed structure that is of particular interest to
us: the separated product and separated function space. Given two
nominal sets $X,Y$, their \emph{separated product} is defined as
$X\ast Y = \set{(x,y)\mid x\in X,\;y\in Y,\; x\# y}$. One can show
that the separated product defines a monoidal structure on
$\Nom$. Moreover, it has a right adjoint $\magicwand$, called the
\emph{separated function space} or \emph{separated closed structure}.
Its elementwise definition is a bit complicated, but we will not need
it here. See {\cite[Def.~3.2]{Clouston2013}} for the details.

It is worth mentioning that a special case of separated closure is
\emph{atom abstraction}: For every nominal set $X$, there is a nominal
set $\aA\magicwand X$, sometimes also written $[\aA]X$, whose elements
are pairs $(a,x)$ up to $\alpha$-equivalence: here $(a,x)$ is
equivalent to $(b,y)$ if for fresh $c$, $(a\,c)\bul x = (b\,c)\bul y$.
An equivalence class of such pairs is usually written $a.x$ or
$\langle a\rangle x$, and represents a \emph{binder}. This is what is
used to model ``syntax with binders'', and is a large part of what
makes the category $\Nom$ interesting.

\section{Extended $\Pi$ action}
\label{sec:large-action}

Recall that $\PI$ is the group of finitely supported permutations of
$\aA$. Let $\PIL$ be the group of all permutations of $\aA$ (not
necessarily finitely supported). If $X$ is a nominal set, then the
group $\PI$ acts on $X$ by definition. Perhaps surprisingly, there is
also a canonical action of $\PIL$ on $X$ \cite{Gadducci2006,milius2016}. The following lemma is
instrumental in proving this.

\begin{xlem}\label{LemSuppAct}
  Let $X$ be a nominal set and let $x\in X$. Let $\pi' ,\pi'' \in \Pi$
  be two finitely supported permutations. If
  $\pi'\restr{\supp(x)}=\pi''\restr{\supp(x)}$ then $\pi' \bul x = \pi''
  \bul x$.
\end{xlem}

\begin{proof}
  Recall that if $\supp(\pi)$ is disjoint from $\supp(x)$, then
  $\pi\bul x = x$.  Now assume
  $\pi'\restr{\supp(x)}=\pi''\restr{\supp(x)}$ and consider
  $\pi=\left( \pi''\right)\inv \pi'\in \Pi$. Then
  $\pi\restr{\supp(x)}=\id\restr{\supp(x)}$, since for all $a\in
  \supp(x)$, we have $\pi'(a)=\pi''(a)$. So we have $\supp(\pi)\cap
  \supp(x)=\emptyset$, which implies $\pi\bul x=x$, and hence
  $\pi'\bul x=\pi''\bul x$, as claimed.
\end{proof}

We can now show that the $\PI$-action on a nominal set can be
extended to $\PIL$.

\begin{xprop}\label{prop:pi-extension}
  Let $(X,\bul)$ be a nominal set. Then the action
  $\bul:\PI\times X \to X$ can be canonically extended to
  an action $\bulL:\PIL\times X \to X$.
\end{xprop} 

\begin{proof}
  See Appendix~\ref{app:proofs}.
\end{proof}

We note that if $\sigma$ happens to be finitely supported, then
$\sigma\bulL x$ and $\sigma\bul x$ agree.

\section{Permutation Induced Functors}

\subsection{The functor $\overline{\iota}$ induced by injective atom renaming}
\label{ssec:overline-iota}

In this section, we consider what happens when we move from one atom
set to another. So let $\aA$ and $\bB$ be two atom sets and consider
an injective function $\iota:\aA\to\bB$. We first note that this induces
a group homomorphism from $\PIA$ to $\PIB$.

\begin{xdefinition}\label{def:PIA-iota}
  Let $\aA$, $\bB$ be atom sets and let $\iota:\aA\to\bB$ be an injective
  function. We define $\PI_\iota:\PIA \to \PIB$ to be the unique
  group homomorphism such that $\PI_\iota \left( (a\,a') \right)=(b\,b')$
  whenever $b=\iota(a)$ and $b'=\iota(a')$. Equivalently, for any permutation
  $\pi\in \PIA$ and any atom $b$, we have
  \[
  \PI_\iota(\pi) (b)=
  \begin{cases} \iota\left( \pi\left( \iota\inv (b) \right)\right) &
    \text{if } b\in \img(\iota)\\ b & \text{otherwise.}
  \end{cases}
  \]
\end{xdefinition}

The above group homomorphism allows us to turn any $\PIB$-action
on a set $X$ into a $\PIA$-action on $X$ by $\pi\diamond x=
\PI_\iota(\pi)\bul x$. Indeed, this defines a functor
$\overline{\iota}:\NomB\to\NomA$.

\begin{xdefinition}
  Given atom sets $\aA$, $\bB$ and an injective function
  $\iota:\aA\to\bB$, the functor $\overline{\iota}:\NomB\to\NomA$ is defined
  as follows.

\begin{itemize}
\item{} On objects: For $(X,\bul)\in\NomB$, define
  $\overline{\iota}(X,\bul)=(X,\diamond)$ where $\pi\diamond x=
  \PI_\iota(\pi)\bul x$ for all $\pi\in \PIA$.
\item{} On morphisms: Given an equivariant map $f:(X,\bul) \to
  (Y,\bul)$ in $\NomB$, we define $\overline{\iota}(f):(X,\diamond)\to
  (Y,\diamond)$ in $\NomA$ by $\overline{\iota}(f)(x)=f(x)$. In other
  words, $\overline{\iota}(f)$ and $f$ have the same underlying function
  $X\to Y$.
\end{itemize}

To show that this is a well-defined functor, we first need to show
that $(X,\diamond)$ is a nominal set. First, it is clear that
$\diamond$ is an action since $\PI_{\iota}$ is a group homomorphism. We
must show that every $x\in X$ is finitely supported. We claim that
$S=\set{a\in \aA \mid \iota(a)\in \supp_{\bB}(x)}$ is a support of $x$.
Indeed, for $a',a''\not\in S$, we have $(a'a'')\diamond x =
(\iota(a'),\iota(a''))\bul x =x$, since $\iota(a'),\iota(a'')\not\in \supp_\bB(x)$.
Therefore, $x$ is finitely supported, and $(X,\diamond)$ is a
well-defined nominal set.  Note that we have proved that
$\supp_{\aA}(x)\subseteq S$, but it is not hard to show that actually
$\supp_{\aA}(x)=S$.

Next, we need to show that $\overline{\iota}$ is well-defined on
morphisms. So given $f:(X,\bul) \to (Y,\bul)$ in $\NomB$, we
must show that $\overline{\iota}(f)$ is equivariant. Therefore, consider
any $\pi \in \PIA$. We must show that
$\overline{\iota}(f)(\pi\diamond x) = \pi \diamond \overline{\iota}(f)(x)$.
Indeed, we have:
\[
\pi \diamond \overline{\iota}(f)(x) = \pi \diamond f(x)
 =  \PI_\iota(\pi) \bul f(x)
 =  f(\PI_\iota(\pi) \bul x)
 =  f(\pi \diamond x)
 =  \overline{\iota} (f)(\pi \diamond x)
\]
So, $\overline{\iota}(f)$ is equivariant. The functoriality of
$\overline{\iota}$ is then trivial, since the composition in both
categories is just composition of the underlying functions, which are
the same.
\end{xdefinition}

\subsection{The functor \texorpdfstring{$\muhat$}{μ} induced by bijective atom renaming}
\label{ssec:muhat-bijective}

An important special case of Section~\ref{ssec:overline-iota} arises when
the atom renaming function is a bijection $\mu:\aA\to\bB$. In this
case, $\PI_{\mu}:\PIA\to\PIB$ is an isomorphism of
groups, and $\overline{\mu}:\NomB\to\NomA$ is an isomorphism of
categories. Moreover, we have $\PI_{\mu}(\pi) = \mu\pi
\mu\inv$ and $\pi\diamond x = (\mu\pi \mu\inv)\bul x$,
for all $\pi\in\PIA$.

The operation that maps $\mu:\aA\to\bB$ to
$\overline{\mu}:\NomB\to\NomA$ is contravariant. Since $\mu$ is
invertible, it is often more convenient to work with its inverse. We
therefore define
\[
\muhat = \overline{\mu\inv}:\NomA\to\NomB.
\]
Since we use the functor $\muhat:\NomA\to\NomB$ a lot, it is useful to
spell out its definition for the record:

\begin{xdefinition}\label{def:muhat}
  For a bijection $\mu:\aA\to\bB$, the functor $\muhat:\NomA\to\NomB$
  is given as follows:
  \begin{itemize}
  \item On objects: For every $(X,\bul)\in\NomA$, we have
    $\muhat(X,\bul)=(X,\diamond)$, where
    \[
    \pi\diamond x = (\mu\inv \pi \mu) \bul x.
    \]
    
  \item On morphisms: For every equivariant function $f:(X,\bul)\to
    (Y,\bul)$, the morphism $\muhat(f):(X,\diamond)\to (Y,\diamond)$
    is defined by $\muhat(f)=f$. 
  \end{itemize}
\end{xdefinition}

\begin{xremark}
  The fact that $\muhat$ acts as the identity on the underlying
  functions of the morphisms was already proved in
  Section~\ref{ssec:overline-iota}, but it is even easier to see in the
  present context, because we have
  \[
  f(\pi\diamond x) =f\big( (\mu\inv \pi \mu) \bul
  x\big) = (\mu\inv \pi \mu) \bul f(x)=\pi\diamond
  f(x).
  \]
  Therefore, $f$ is equivariant with respect to both the
  $\bul$-action and the $\diamond$-action.
\end{xremark}

Recall from Section~\ref{sec:large-action} that every bijection
$\sigma:\aA\to\aA$ (whether or not it is finitely supported) acts on
every $\aA$-nominal set via $\bulL$. One may ask whether the action
$\sigma\bulL(-):X\to X$ is a natural transformation. However, this is
not the case. In fact, $\sigma\bulL(-):X\to X$ is not even a morphism
of $\NomA$, since it is not equivariant. Equivariance would require
that $\sigma\bulL(\pi\bul x) = \pi\bul(\sigma\bulL x)$, which in
general only holds if $\sigma$ and $\pi$ commute. We do, however, have
the following:

\begin{xprop}\label{DotBlankMapNatIsoProp}
  Let $\sigma:\aA\to\aA$ be a bijection. Then $\DBL{\sigma}:\sigmahat(X)\to X$
  is a natural isomorphism.
\end{xprop}

\begin{proof}
  See Appendix~\ref{app:proofs}.
\end{proof}

The next lemma deals with the situation where we have two permutations
$\sigma_1,\sigma_2:\aA\to\aA$. 

\begin{xlem}\label{lem:mu1mu2}
 Consider permutations $\sigma_1,\sigma_2:\aA\to\aA$, and let
  $\sigma=\sigma_2\circ\sigma_1$.   Let $(X,\bul)$ be an object of $\NomA$, and let
  $\sigmahat_1(X,\bul)=(X,\diamond)$ as in Definition~\ref{def:muhat}. Then the following diagram commutes:
  \[
  \begin{tikzpicture}[scale=0.7]
    \node (1) at (0,3) {$\sigmahat_2(\sigmahat_1(X))$};
    \node (2) at (-3,0) {$\sigmahat_2(X)$};
    \node (3) at (0,1.3) {$\sigmahat(X)$};
    \node (4) at (3,0) {$\sigmahat_1(X)$};
    \node (5) at (0,-3) {$X$};

    \path[every node/.style={font=\sffamily\small}]
    (1) [black,thick] edge [double distance=2pt] (3)
    (1) [-stealth,black,thick] edge node [left] {$\sigmahat_2(\DBL{\sigma_1})$} (2)
    (1)  edge node [right] {$\DDL{\sigma_2}$} (4)
    (2)  edge node [left] {$\DBL{\sigma_2}$} (5)
    (3)  edge node [right,near start] {$\DBL{\sigma}$} (5)
    (4)  edge node [right] {$\DBL{\sigma_1}$} (5)
    ;
  \end{tikzpicture}
  \]
\end{xlem}

\begin{proof}
  The morphisms along the left and right of this diagram map $x$ to
  $\sigma_2\bulL(\sigma_1\bulL x)$ and $\sigma_1\bulL(\sigma_2\diamondL
  x)$, respectively. By the definition of the various actions, both
  are equal to $\sigma\bulL x$. Alternatively, we could have also
  observed that the outer square commutes by naturality.
\end{proof}

\begin{xlem}\label{lem:permutationconjugatemap}
  Let $\mu:\aA \to \bB$ and $\sigma:\aA \to \aA$ be bijections.  Let
  $(X,\bul)$ be an object of $\NomA$, and let
  $\muhat(X,\bul)=(X,\diamond)$ as in Definition~\ref{def:muhat}.
  Then $\DDL{\mu\sigma\mu\inv} = \widehat{\mu}(\DBL{\sigma})$, i.e., the
  following diagram commutes in $\NomB$.
  \[
  \begin{tikzpicture}[scale=0.7]
    \node (1) at (0,3) {$\widehat{\mu\sigma}(X)$};
    \node (2) at (-3,0) {$\widehat{\mu}(\sigmahat(X))$};
    \node (3) at (3,0) {$\widehat{\mu\sigma\mu\inv}(\widehat{\mu}(X))$};
    \node (4) at (0,-3) {$\muhat(X)$};

    \path[every node/.style={font=\sffamily\small}]
    (1) [black,thick] edge [double distance=2pt] (2)
    (1) [black,thick] edge [double distance=2pt] (3)
    (2) [-stealth,black,thick] edge node [left] {$\muhat(\DBL{\sigma})$} (4)
    (3)  edge node [right] {$\DDL{\mu\sigma\mu\inv}$} (4)
    ;
  \end{tikzpicture}
  \]
\end{xlem}

\begin{proof}
  Let $x\in \widehat{\mu\sigma}(X)$. Then $(\mu\sigma\mu\inv)\diamondL x
   = \mu\inv (\mu\sigma\mu\inv) \mu \bul x = \sigma \bul x$. Also, since the 
   $\muhat$ functor acts trivially on the underlying function, 
   $\muhat(\DBL{\sigma})$ acts the same as $\DBL{\sigma}$ on elements.
\end{proof}

In this paper, we will rarely deal with two completely arbitrary atom
sets $\aA$ and $\bB$; instead, we will most often use the functor
$\sigmahat$ in the context of subsets of a fixed atom set $\aA$.

\begin{xnotation}
  Fix an atom set $\aA$. Whenever $S\subseteq\aA$ is a finite subset, we
  write $\aA_S=\aA\setminus S$.
\end{xnotation}

Note that every permutation $\sigma:\aA\to\aA$
(whether or not it is finitely supported) restricts to a bijective
function $\sigma: \aA_S \xrightarrow{\cong} \aA_{\sigma\bulL S}$. By
Definition~\ref{def:muhat}, this in turn induces an isomorphism of categories
\[
\sigmahat :\NomAX{S}\xrightarrow{\cong}\NomAX{\sigma\bulL S}.
\]

\begin{xremark}
  In this section, we considered that a bijection of atom sets
  $\mu:\aA\to\bB$ induces an isomorphism of categories
  $\muhat:\NomA\to\NomB$. In particular,
the categories $\NomA$ and $\NomB$ for two different countable atom sets are always equivalent. It is an interesting observation that even if we were to allow atom sets of arbitrary infinite cardinalities, the
categories $\NomA$ and $\NomB$ still end up being equivalent. We do not need this more general fact in this paper, and we will not prove it here. However, it justifies that working with a countable atom set is
without loss of generality.
\end{xremark}

\section{Nominal families}

Fix a set $\aA$ of atoms, and let $X$ be an object in $\NomA$.
Consider the slice category $\NomA/X$ and an object in it:
\[
\begin{tikzpicture}
  \node (1) at (0,0) {$Y$};
  \node (2) at (0,-1.5) {$X$};
  
  \path[every node/.style={font=\sffamily\small}]
  (1) [-stealth,black,thick] edge node [midway,right] {$p$} (2);
\end{tikzpicture}
\]
Let $x\in X$ and let $Y_x$ denote the fiber over $x$, i.e.,
$Y_x=\set{y\in Y \mid p(y)=x}$. Note that $Y_x$ is not in general an
object of $\NomA$, because the $\pi$-action is not well-defined on it:
For $y\in Y_x$, we have $\pi\bul y\in Y_{\pi\bul x}$ instead of
$\pi\bul y\in Y_x$.

However, $Y_x$ is an object of $\NomAX{x}$, where $\aA_x=\aA\setminus
\supp(x)$. Specifically, let $\iota:\aA_x\hookrightarrow\aA$ be the subset
inclusion, and recall from Definition~\ref{def:PIA-iota} that every
permutation $\pi\in\PI_{\aA_x}$ can be canonically extended to a
permutation $\PI_\iota(\pi)\in\PIA$. Then the action of
$\pi\in\PI_{\aA_x}$ on an element $y\in Y_x$ is defined as
\[
\pi \bul_{Y_x} y = \PI_\iota(\pi)\bul_Y y.
\]
We call $\NomAX{x}$ the category of \textit{nominal sets relative to
  $x$}.

Conversely, consider a family $(Y_x)_{x\in X}$ of nominal sets,
where each $Y_x$ is an object of $\NomAX{x}$. With what additional
structure must we equip this family, so that it will give
rise to an object $p:Y\to X$ of the slice category? Clearly, we want
$Y=\sum_{x\in X} Y_x$, where $\sum$ denotes the disjoint union.
The additional structure we need on the family $(Y_x)_{x\in X}$ is a
``global'' action $\bulG$ such that for all $\pi\in\PIA$, we have
$\pi\bulG(-) : Y_x \to Y_{\pi\bul x}$, subject to some
axioms. This motivates the following definition.

\begin{xdefinition}\label{def:nominal-family}
  Let $X$ be a nominal set over the atom set $\aA$. A \emph{nominal
  family over $X$} is a family $(Y_x)_{x\in X}$ such that each $Y_x$
  is a nominal set in $\NomAX{x}$, together with a family of
  operations $\pi\bulG(-) : Y_x \to Y_{\pi\bul x}$, satisfying
  the following properties:
  \begin{enumerate}\alphalabels
  \item\label{b-b} $\pi_1\pi_2\bulG y =
    \pi_1\bulG(\pi_2\bulG y)$.
  \item\label{b-c} If $\pi\# x$, then for all $y\in Y_x$, we have
    $\pi\bulG y=\pi\bul y$.
  \end{enumerate}
\end{xdefinition}

Property (\ref{b-b}) is like the usual property of an action, and
property (\ref{b-c}) states that the ``global'' action $\bulG$
coincides with the ``local'' action $\bul$ on the nominal set $Y_x$
when both are defined. Note that we did not require $\id\bulG y =
y$, because this follows from (\ref{b-c}) by setting $\pi=\id$.

Given two nominal families $(Y_x)_{x\in X}$ and $(Z_x)_{x\in X}$,
consider a family of morphisms $f_x:Y_x\to Z_x$ where each $f_x$ is a
morphism in $\NomAX{x}$. We say that such a family of morphisms is
\emph{equivariant} if it respects the global action, i.e., if for all
$x\in X$, $y\in Y_x$, and $\pi\in\PIA$, we have
\[
\pi\bulG (f_x(y)) = f_{\pi\bul x}(\pi\bulG y).
\]
The \emph{category of nominal families over $X$} has nominal families
as the objects and equivalent families of morphisms as the maps. We
denote it by $\Nom^X$. The following proposition is easy to prove.

\begin{xprop}\label{prop:nom-x-equivalent}
  Let $X\in\NomA$. The categories $\Nom^X$ and $\Nom/X$ are equivalent.
\end{xprop}

\begin{proof}
  We already saw that an object $p:Y\to X$ of the slice category gives
  rise to a nominal family $(Y_x)_{x\in X}$. Conversely, if
  $(Y_x)_{x\in X}$ is a nominal family, let $Y=\sum_{x\in X} Y_x =
  \set{(x,y)\mid x\in X, y\in Y_x}$ be the disjoint union with $p:Y\to
  X$ defined by $p(x,y) = x$. We can equip $Y$ with the action
  $\pi\bul(x,y) = (\pi\bul x, \pi\bulG y)$. Every element
  $(x,y)$ is finitely supported by
  $\supp_{\aA}(x)\cup\supp_{\aA_{x}}(y)$.  It is routine to verify
  that it is a well-defined object of $\Nom/X$, and that the two
  constructions respect morphisms and are mutually inverse up to
  isomorphism.
\end{proof}

Consider a nominal family $(Y_x)_{x\in X}$, some permutation
$\pi\in\PIA$ and some $x\in X$. We may ask what is the
relationship between the fibers $Y_x$ and $Y_{\pi\bul x}$. The
global action gives us functions $\pi\bulG(-):Y_x\to Y_{\pi\bul
  x}$ and $\pi\inv\bulG(-):Y_{\pi\bul x}\to Y_x$. Since these
functions are each other's inverses, it follows that $Y_x$ and
$Y_{\pi\bul x}$ are isomorphic as sets. Are they isomorphic as
nominal sets?  The question does not really make sense, because $Y_x$
and $Y_{\pi\bul x}$ are objects in two different categories, namely
in $\NomAX{x}$ and $\NomAX{\pi\bul x}$, respectively. On the other
hand, as we saw in Section~\ref{ssec:muhat-bijective}, there is an
isomorphism of categories $\pihat:\NomAX{x}\to\NomAX{\pi\bul
  x}$. Modulo this isomorphism, the objects $Y_x$ and $Y_{\pi\bul
  x}$ are isomorphic:

\begin{xprop}\label{PiEquivIsoLem}
  Given a nominal family $(Y_x)_{x\in X}$, some $\pi\in\PIA$ and
  $x\in X$. Then the function $\eqfam{\pi}{x} : \pihat(Y_x)\to
  Y_{\pi\bul x}$ defined by $\eqfam{\pi}{x}(y) = \pi\bulG y$ is an
  equivariant isomorphism.
\end{xprop}

\begin{proof}
  Clearly $\eqfam{\pi}{x}$ is invertible, with its inverse given by
  the action of $\pi\inv$. We must show that $\eqfam{\pi}{x}$ is
  equivariant. Recall from Section~\ref{ssec:muhat-bijective} that
  $\pihat(Y_x)=(Y_x,\diamond)$, with the action $\rho\diamond y =
  (\pi\inv\rho\pi)\bul y$.  Consider some $y\in Y_x$ and
  $\rho\in\PI_{\pi\bul x}$. We have
  \[
  \eqfam{\pi}{x}(\rho\diamond y)
  = \pi\bulG((\pi\inv\rho\pi)\bul y)
  = \pi\bulG(\pi\inv\rho\pi)\bulG y
  = \pi\bulG\pi\inv\bulG\rho\bulG\pi\bulG y
  = \rho\bulG\pi\bulG y
  = \rho\bul(\pi\bulG y)
  = \rho\bul \eqfam{\pi}{x}(y).
  \]
  Here, in the relevant steps, we were able to replace $\bul$ by
  $\bulG$ and vice versa because $(\pi\inv\rho\pi)\# x$ and
  $\rho\#(\pi\bul x)$.
\end{proof}

\begin{xremark}\label{rk:NomFamPhiMap}
  We can equivalently define a nominal family in terms of a family of
  maps $\eqfam{\pi}{x}$, rather than a global action $\bulG$. In
  this case, a nominal family over $X$ consists of a family
  $\Fam{Y_x}_{x\in X}$, where $Y_x\in\NomAX{x}$, together with an
  equivariant map $\eqfam{\pi}{x}:\pihat(Y_x)\to Y_{\pi\bul x}$ for
  each $\pi\in \PIA$ and $x\in X$, satisfying the following
  conditions:
  \begin{enumerate}\alphalabels
  \item\label{ax-c} $\eqfam{\pi_1\pi_2}{x}=\eqfam{\pi_1}{\pi_2\bul x}\circ \widehat{\pi_1}(\eqfam{\pi_2}{x})$.
  \item\label{ax-d} If $\pi\# x$, then $\eqfam{\pi}{x}(y) = \pi
    \bul y$ for all $y\in \pihat(Y_x)$.
  \end{enumerate}
  The equivalence between this definition and
  Definition~\ref{def:nominal-family} is straightforward, as both sets
  of conditions refer to the same sets and elements. Also, due to
  Proposition~\ref{PiEquivIsoLem}, conditions
  (\ref{ax-c}) and (\ref{ax-d}) automatically imply that the maps
  $\eqfam{\pi}{x}$ are equivariant and invertible.

  In this case, a family of morphisms $f_x:Y_x\to Z_x$
  is equivariant if and only if the following diagram commutes for all
  $x\in X$ and $\pi\in\PIA$:
  \[
  \begin{tikzpicture}
    \node (ltop) at (0,0) {$\pihat(Y_x)$};
    \node (rtop) at (3,0) {$\pihat(Z_x)$};
    \node (lbot) at (0,-1.5) {$Y_{\pi\bul x}$}; 
    \node (rbot) at (3,-1.5) {$Z_{\pi\bul x}$}; 
    
    \draw [-stealth][black,thick] (ltop) -- (rtop) node[midway,above] {$\pihat(f_x)$};
    \draw [-stealth][black,thick] (ltop) -- (lbot) node[midway,left] {$\eqfams{Y}{\pi}{x}$};
    \draw [-stealth][black,thick] (lbot) -- (rbot) node[midway,above] {$f_{\pi\bul x}$};
    \draw [-stealth][black,thick] (rtop) -- (rbot) node[midway,right] {$\eqfams{Z}{\pi}{x}$};
  \end{tikzpicture}
  \]
\end{xremark}

\section{Monoidal structure}
\label{sec:monoidal}

In this section, we will discuss the relation between the monoidal
structure on $\Nom$ and the monoidal structure on $\Nom/X$.

\begin{xtheorem}\label{thm:monoidal-nom-slash-x}
    Every monoidal structure $\otimes$ on $\Nom$ induces a
    monoidal structure $\otimes_X$ on $\Nom/X$, which is uniquely determined up to coherent isomorphism.
\end{xtheorem}

Before we prove the theorem, we need to review some facts about how
monoidal structures transport along isomorphisms.

\subsection{Monoidal structures and isomorphisms of categories}
\label{ssec:monoidal-iso-cats}

Recall that a functor $F:\Cc\to\Dd$ between monoidal categories is
said to be a {\em (strong) monoidal functor} if it is equipped with
natural isomorphisms $\mm^0 : FI\to I$ and $\mm^{A,B}:F(A\otimes B)\to
FA\otimes FB$ that are compatible with the associators $\alpha_{A,B,C}$ and left and right unitors
$\lambda_A$ and $\rho_A$ in
the obvious way. Note that in the literature, it is more common to
define the arrows $\mm^0$ and $\mm^{A,B}$ to point in the opposite
direction, but since they are isomorphisms, this does not matter. It
is more convenient for our purposes to define them as above. Also,
since we are never interested in lax (nor op-lax) monoidal functors,
we will drop the adjective ``strong''.

Consider a category $\Cc$ that is equipped with two monoidal
structures $(\otimes,I,\lambda,\rho,\alpha)$ and
$(\otimes',I',\lambda',\rho',\alpha')$. In this case, an
\emph{isomorphism} between the monoidal structures is a monoidal
functor $F:(\Cc,\otimes,I,\lambda,\rho,\alpha) \to
(\Cc,\otimes',I',\lambda',\rho',\alpha')$ whose underlying functor is
the identity. Specifically, this means that there are natural
isomorphisms $\mm^0 : I\to I'$ and $\mm^{A,B}:A\otimes B\to
A\otimes' B$ that are compatible with the associators and unitors.

We now consider how monoidal structures are transported along
isomorphisms of categories. Most of what we say below could be
generalized to equivalences of categories, but for the purposes of
this paper, it is sufficient to consider isomorphisms.

Let $F:\Cc\to\Dd$ be an isomorphism of categories. Then every monoidal
structure $(\otimes,I,\lambda,\rho,\alpha)$ on $\Cc$ canonically
induces a monoidal structure
$(\otimes_F,I_F,\lambda_F,\rho_F,\alpha_F)$ on $\Dd$ in the obvious
way. For example, $A\otimes_F B = F(F\inv(A)\otimes F\inv(B))$, $I_F =
FI$, etc.

In general, if $F,G:\Cc\to\Dd$ are two isomorphisms of categories,
then the two induced monoidal structures $\otimes_F$ and $\otimes_G$
may not be isomorphic. For example, consider $\Cc=\Dd=\Set\times\Set$ with
the monoidal structure $(A,B)\otimes(A',B')=(A\times A',B+B')$. Then
the functors $F(A,B)=(A,B)$ and $G(A,B)=(B,A)$ induce two
non-isomorphic monoidal structures on $\Set\times\Set$.

If, on the other hand, the functors $F,G$ are naturally isomorphic to
each other, say via $\eta^A:FA\to GA$, then the monoidal structures
induced on $\Dd$ by $F$ and $G$ are isomorphic in a canonical way.

\subsection{Transport of monoidal structure from $\NomA$ to $\NomB$}

We now want to show that if $\aA$ and $\bB$ are two (countably
infinite) atom sets, then any monoidal structure on $\Nom_{\aA}$
determines a monoidal structure on $\Nom_{\bB}$, and moreover, the
latter structure is uniquely determined up to isomorphisms of
monoidal structures (and in particular, it does not depend on any
particular choice of a bijection between $\aA$ and $\bB$).

So consider a monoidal structure on $\NomA$. Clearly, any bijection
$\gamma:\aA\to\bB$ induces an isomorphism between $\NomA$ and $\NomB$,
and therefore induces a monoidal structure on $\NomB$. What we must
show is that the latter monoidal structure is independent of $\gamma$.

\begin{xprop}\label{ISOProp}
Let $\gamma_1,\gamma_2:\aA\to \bB$ be two bijections and let
$\otimes_1$ and $\otimes_2$ be the monoidal structures on $\NomB$
induced by $\gamma_1$ and $\gamma_2$, respectively.  The monoidal structures $\otimes_1$ and $\otimes_2$ are isomorphic.
\end{xprop}

\begin{proof}
  Define $\sigma:\bB \to \bB$ by $\sigma=\gamma_2\circ (\gamma_1)\inv$.
  Consider $\sigmahat:\NomB\to\NomB$ as in
  Section~\ref{ssec:muhat-bijective}. Note that $\sigmahat$ is naturally
  isomorphic to the identity functor $\id:\NomB\to\NomB$ by
  Proposition~\ref{DotBlankMapNatIsoProp}. Therefore, by whiskering,
  $\sigmahat\circ\widehat\gamma_1=\widehat\gamma_2$ is naturally isomorphic
  to $\id\circ\widehat\gamma_1=\widehat\gamma_1$. As noted in
  Section~\ref{ssec:monoidal-iso-cats}, this implies that the monoidal
  structures on $\NomB$ induced by $\widehat\gamma_2$ and
  $\widehat\gamma_1$ are isomorphic, as claimed.
\end{proof}

\begin{xremark}\label{rem:EtaMap}
  We write $\eta^{A,B}_{\gamma_1,\gamma_2} : A\otimes_1 B\to A\otimes_2 B$
  for the natural isomorphism described in
  Proposition~\ref{ISOProp}. It satisfies the following diagram, where
  $\sigma=\gamma_2\circ (\gamma_1)\inv$:
  \[
  \begin{tikzpicture}
    \node (ltop) at (0,0) {$\widehat{\sigma}(A\otimes_1 B)$};
    \node (rtop) at (3,0) {$\widehat{\sigma}(A)\otimes_2 \widehat{\sigma}(B)$};
    \node (lbot) at (0,-1.5) {$A\otimes_1 B$}; 
    \node (rbot) at (3,-1.5) {$A\otimes_2 B$.}; 
    
    \draw [black,  thick, double distance = 1.5pt] (ltop) -- (rtop) node[midway,above] {};
    \draw [-stealth][black,  thick] (ltop) -- (lbot) node[midway,left] {$\DBL{\sigma}$}node[midway,right] {$\cong$};
    \draw [-stealth][black,  thick] (lbot) -- (rbot) node[midway,above] {$\eta_{\gamma_1,\gamma_2}^{A,B}$};
    \draw [-stealth][black,  thick] (rtop) -- (rbot) node[midway,left] {$\cong$} node[midway,right] {$\DBL{\sigma}\otimes_2 \DBL{\sigma}$};
  \end{tikzpicture}
  \]
  Similarly, we write $\eta^0_{\gamma_1,\gamma_2} : I_1\to I_2$ for the
  corresponding isomorphism between the units of the two monoidal
  structures. It satisfies the following:
  \[
  \begin{tikzpicture}
    \node (ltop) at (0,0) {$\widehat{\sigma}(I_1)$};
    \node (rtop) at (3,0) {$I_2$};
    \node (lbot) at (0,-1.5) {$I_1$}; 
    \node (rbot) at (3,-1.5) {$I_2$.}; 
    
    \draw [black,  thick, double distance = 1.5pt] (ltop) -- (rtop);
    \draw [-stealth][black,  thick] (ltop) -- (lbot) node[midway,left] {$\DBL{\sigma}$}node[midway,right] {$\cong$};
    \draw [-stealth][black,  thick] (lbot) -- (rbot) node[midway,above] {$\eta_{\gamma_1,\gamma_2}^{0}$};
    \draw [black,  thick, double distance = 1.5pt] (rtop) -- (rbot);
  \end{tikzpicture}
  \]
\end{xremark}

Proposition~\ref{ISOProp} ensures that the monoidal structures induced
on $\NomB$ by various different bijections $\aA\to\bB$ are all
isomorphic to each other. We will also need these isomorphisms to be
coherent, i.e., any three such monoidal structures are isomorphic in a
unique way. The following lemma shows this.

\begin{xlem}\label{lem:EtaCompose}
  Let $\gamma_1,\gamma_2,\gamma_3:\aA\to \bB$ be three bijections and let
  $\otimes_1$, $\otimes_2$ and $\otimes_3$ be the monoidal structures
  on $\NomB$ induced by $\gamma_1$, $\gamma_2$, and $\gamma_3$,
  respectively.  Then $\eta^{A,B}_{\gamma_1,\gamma_3} =
  \eta^{A,B}_{\gamma_2, \gamma_3} \circ \eta^{A,B}_{\gamma_1,\gamma_2}$.
\end{xlem}

\begin{proof}
  See Appendix~\ref{app:proofs}.
\end{proof}

\begin{xlem}\label{lem:eta-pihat}
  Let $\gamma_1,\gamma_2:\aA\to\bB$ be bijections, inducing monoidal
  structures $\otimes_1$ and $\otimes_2$ on $\NomB$. Let
  $\pi:\bB\to\cC$ be another bijection, and let $\otimes'_1$, and
  $\otimes'_2$ be the monoidal structures induced on $\NomC$ by
  $\pi\gamma_1$, and $\pi\gamma_2$, respectively.  
  Then the map $\eta^{A,B}_{\gamma_1,\gamma_2}$ is compatible with
  the $\widehat{\pi}$ functor, i.e., the following diagram commutes.
  \[
  \begin{tikzpicture}
    \node (1) at (0,0) {$\widehat{\pi}(A\otimes_1 B)$};
    \node (2) at (5,0) {$\widehat{\pi}(A\otimes_2 B)$};
    \node (3) at (0,-1.5) {$\widehat{\pi}(A)\otimes'_1 \widehat{\pi}(B)$};
    \node (4) at (5,-1.5) {$\widehat{\pi}(A)\otimes'_2 \widehat{\pi}(B)$};
    
    \path[every node/.style={font=\sffamily\small}]
    (1) [black,thick] edge [double distance=2pt] (3)
    (2) [black,thick] edge [double distance=2pt] (4)
    
    (1) [-stealth,black,thick] edge node [midway,below]{$\widehat{\pi}(\eta_{\gamma_1,\gamma_2}^{A,B})$} (2)
    (3) [-stealth,black,thick] edge node [midway,below]{$\eta_{\pi\gamma_1,\pi\gamma_2}^{A,B}$} (4);
  \end{tikzpicture}
  \]
\end{xlem}

\begin{proof}
  See Appendix~\ref{app:proofs}.
\end{proof}

\subsection{\texorpdfstring{$\muhat$}{μ} is a monoidal functor}
\label{ssec:mmap}

As usual, let $\otimes$ be a monoidal structure on $\NomA$.  Consider
bijections $\delta_1:\aA\to\bB$ and $\delta_2:\aA\to\cC$, and let
$\otimes_1$ and $\otimes_2$ be the monoidal structures induced by
$\delta_1$ and $\delta_2$ on $\NomB$ and $\NomC$, respectively.

Now consider another arbitrary bijection $\mu:\bB\to\cC$, inducing an
isomorphism of categories $\muhat:\NomB\to\NomC$. Note that we do not
require $\delta_2 = \mu\delta_1$. Therefore, by transporting $\otimes_1$
along $\muhat$, we obtain another monoidal structure on $\NomC$, which
we call $\otimes'_2$, and which is not necessarily the same as
$\otimes_2$. By definition, $\muhat$ strictly maps $\otimes_1$ to
$\otimes'_2$, and therefore it certainly is a monoidal functor with
respect to $\otimes_1$ and $\otimes'_2$. Less obvious is the fact that
$\muhat$ is also canonically a monoidal functor with respect to
$\otimes_1$ and $\otimes_2$. The purpose of this subsection is to
spell out how.

The key is Proposition~\ref{ISOProp}, which shows that the two monoidal
structures on $\NomC$, $\otimes_2$ and $\otimes'_2$, while not the
same, are canonically isomorphic. Therefore, $\muhat$ can be
canonically equipped with the structure of a monoidal functor with
respect to $\otimes_1$ and $\otimes_2$. Explicitly, we define a family
of maps $\mm_{\mu}:\muhat(A\otimes_1
B)\to\muhat(A)\otimes_2\muhat(B)$ by the following diagram:
\vspace{-0.25cm}
\[
\begin{tikzpicture}
  \node (1) at (0,0) {$\muhat(A\otimes_1 B)$};
  \node (2) at (4,0) {$\muhat(A)\otimes'_2\muhat(B)$};
  \node (3) at (8,0) {$\muhat(A)\otimes_2\muhat(B)$.};
  
  \path[every node/.style={font=\sffamily\small}]
  (1) [black,thick] edge [double distance=2pt] (2)
  (2) [-stealth,black,thick] edge node [midway,above] {$\eta_{\mu\circ\delta_1,\delta_2}$} (3)
  (1) [-stealth,black,thick] edge [bend right=15] node [midway,below]{$\mm_{\mu}$} (3);
\end{tikzpicture}
\]
\vspace{-1cm}

\noindent
Similarly, we also define
\vspace{-0.25cm}
\[
\begin{tikzpicture}
  \node (1) at (0,0) {$\muhat(I_1)$};
  \node (2) at (4,0) {$I'_2$};
  \node (3) at (8,0) {$I_2$.};
  
  \path[every node/.style={font=\sffamily\small}]
  (1) [black,thick] edge [double distance=2pt] (2)
  (2) [-stealth,black,thick] edge node [midway,above] {$\eta^{0}_{\mu\circ\delta_1,\delta_2}$} (3)
  (1) [-stealth,black,thick] edge [bend right=15] node [midway,below]{$\mm^{0}_{\mu}$} (3);
\end{tikzpicture}
\]
\vspace{-0.75cm}

\noindent
We note that $\mm_{\mu}$ is a natural isomorphism because the equality
$\muhat(A\otimes_1 B) = \muhat(A)\otimes'_2\muhat(B)$ is natural by
definition of $\otimes'_2$, and $\eta_{\mu\circ\delta_1,\delta_2}$
is natural by Remark~\ref{rem:EtaMap} and the proof of
Proposition~\ref{ISOProp}.  Also, $\mm_{\mu}$ is compatible with the
monoidal associators and unitors. Indeed, the compatibility of the
equality $\muhat(A\otimes_1 B) = \muhat(A)\otimes'_2\muhat(B)$ holds
by the definition of the associators and unitors for $\otimes'_2$,
whereas the compatibility for $\eta_{\mu\circ\delta_1,\delta_2}$
holds because it is an isomorphism of monoidal structures by
Proposition~\ref{ISOProp}. Thus, $\mm_{\mu}$ indeed makes $\muhat$ into a
monoidal functor from $(\NomB,\otimes_1)$ to $(\NomC,\otimes_2)$.

The following lemma establishes a useful property of the $\mm$ maps.

\begin{xlem}\label{lem:mmapcompose}
  Let $\delta_1:\aA\to\bB$, $\delta_2:\aA\to\cC$, and
  $\delta_3:\aA\to\dD$ be three bijections, inducing monoidal
  structures $\otimes_1$, $\otimes_2$, and $\otimes_3$ on $\NomB$,
  $\NomC$, and $\NomD$, respectively. Let $\mu_1:\bB\to \cC$ and
  $\mu_2:\cC\to \dD$ be bijections. Then the following diagram
  commutes. 
  \[
  \begin{tikzpicture}
    \node (1) at (0,0) {$\muhat_1\muhat_2(A\otimes_1 B)$};
    \node (2) at (5,0) {$\muhat_1(\muhat_2(A)\otimes_2\muhat_2(B))$};
    \node (3) at (10,0) {$\muhat_1\muhat_2(A)\otimes_3\muhat_1\muhat_2(B)$};
    \node (4) at (0,-1.25) {$\widehat{\mu_1\mu_2}(A\otimes_1 B)$};
    \node (5) at (10,-1.25) {$\widehat{\mu_1\mu_2}(A)\otimes_3\widehat{\mu_1\mu_2}(B)$};

    \path[every node/.style={font=\sffamily\small}]
    (1) [black,thick] edge [double distance=2pt] (4)
    (3) [black,thick] edge [double distance=2pt] (5)
    (1) [black,thick] edge node [midway,above] {$\muhat_1(\mm_{\mu_2})$} (2)
    (2) [-stealth,black,thick] edge node [midway,above] {$\mm_{\mu_1}$} (3)
    (4) [-stealth,black,thick] edge node [midway,below]{$\mm_{\mu_1\mu_2}$} (5);

  \end{tikzpicture}
  \]
\end{xlem}

\begin{proof}
  See Appendix~\ref{app:proofs}.
\end{proof}

A special case arises when $\delta_1=\delta_2$. In this case, we have
the following:

\begin{xlem}\label{lem:m-pi}
  Let $\delta_1=\delta_2:\aA\to\bB$, so that the induced tensors
  $\otimes_1$ and $\otimes_2$ on $\NomA$ coincide. Consider any
  $\mu:\bB\to\bB$. Then we have
  \[
  \begin{tikzpicture}
    \node (1) at (0,0) {$\muhat(A\otimes_1 B)$};
    \node (3) at (5,0) {$\muhat(A)\otimes_2\muhat(B)$};
    \node (4) at (0,-1.5) {$A\otimes_1 B$};
    \node (6) at (5,-1.5) {$A\otimes_2 B$.};

    \path[every node/.style={font=\sffamily\small}]
    (4) [black,thick] edge [double distance=2pt] (6)
    (1) [-stealth,black,thick] edge node [midway,above] {$\mm_{\mu}$} (3)
    (1) [black,thick] edge node [midway,left]{$\DBL{\mu}$} (4)
    (3) [black,thick] edge node [midway,right] {$\DBL{\mu}\otimes_2\DBL{\mu}$} (6)
    ;
  \end{tikzpicture}
  \]
\end{xlem}

\begin{proof}
  See Appendix~\ref{app:proofs}.
\end{proof}

\subsection{Proof of Theorem~\ref{thm:monoidal-nom-slash-x}}

Consider a monoidal structure $\otimes$ on $\Nom$. We wish to define a
monoidal structure $\otimes_X$ on $\Nom/X$. From
Proposition~\ref{prop:nom-x-equivalent}, we know that $\Nom/X$ is
equivalent to the category $\Nom^X$ of nominal families over $X$.
It therefore suffices to define a monoidal structure on the latter.

\subsubsection{Definition of $A\otimes_X B$}\label{sssec:def-a-times-b}

Consider two nominal families $A = (A_x)_{x\in X}$ and $B =
(B_x)_{x\in X}$. We write $\pi\bulG(-)$ for the ``global'' action
on both of these families. We wish to define a nominal family
$A\otimes_X B = (C_x)_{x\in X}$.

We start by choosing a fixed but arbitrary bijection $\delta_x :
\aA\to\aA_{x}$ for every $x\in X$. This induces an isomorphism of
categories $\widehat\delta_x:\NomA\to\NomAX{x}$, and therefore a
monoidal structure in $\NomAX{x}$. We write $\otimes_x$ for this
monoidal structure (and also $I_x$, $\alpha_x$, etc., as necessary).

Next, we define $C_x = A_x\otimes_x B_x$. To make this into a nominal
family, we must define a global $\pi$-action, i.e., a family of
operations $\pi\bulG(-) : C_x \to C_{\pi\bul x}$. So fix some
permutation $\pi$ and an $x\in X$. For ease of notation, let us write
$y=\pi\bul x$.

Each of the sets $A_x$ and $B_x$ is already equipped with an action,
namely $\pi\bulG(-):A_x\to A_y$ and $\pi\bulG(-):B_x\to B_y$.
By Proposition~\ref{PiEquivIsoLem}, these correspond to equivariant
maps $\eqfams{A}{\pi}{x} : \pihat(A_x)\to A_y$ and $\eqfams{B}{\pi}{x}
: \pihat(B_x)\to B_y$. Since both of these maps are morphisms in the
category $\NomAX{y}$, we may tensor them to get a map
\[
\eqfams{A}{\pi}{x}\otimes_y\eqfams{B}{\pi}{x}
~:~
\pihat(A_x)\otimes_y\pihat(B_x)
~\to~
A_y\otimes_yB_y.
\]
We now define the map $\eqfams{C}{\pi}{x} : \pihat(C_x)\to C_y$ as
follows, using the $\mm$-map of Section~\ref{ssec:mmap}.
\[
  \begin{tikzpicture}[xscale=0.8]
    \node (1) at (-6,0) {$\pihat(C_x)$};
    \node (2) at (-6,-1.25) {$\pihat(A_x\otimes_x B_x)$};
    \node (4) at (0,-1.25) {$\pihat(A_x)\otimes_y\pihat(B_x)$}; 
    \node (5) at (6,-1.25) {$A_y\otimes_yB_y$};
    \node (6) at (6,0) {$C_y$};

    \path[every node/.style={font=\sffamily\small}]
    (1) [black,thick] edge [double distance=2pt] (2)
    (5) [black,thick] edge [double distance=2pt] (6)
    (2) [-stealth,black,thick] edge node [midway,above] {$\mm_{\pi}$} (4)
    (4) [-stealth,black,thick] edge node [midway,above]{$\eqfams{A}{\pi}{x}\otimes_y\eqfams{B}{\pi}{x}$} (5)
    (1) [-stealth,black,thick] edge node [midway,above]{$\eqfams{C}{\pi}{x}$} (6);
\end{tikzpicture}
\]
By Remark~\ref{rk:NomFamPhiMap}, the family of maps $\eqfams{C}{\pi}{x}$ is
sufficient to turn $(C_x)_{x\in X}$ into a nominal family, provided
that it satisfies the remark's properties (\ref{ax-c}) and (\ref{ax-d}).
This is shown in Appendix~\ref{app:prop-a-and-b-1}.
This finishes the proof that $A\otimes_X B$ is a well-defined object.

\subsubsection{Definition of $f\otimes_X g$ and functoriality}

To continue our definition of the monoidal structure on $\Nom/X$, we
must now define the morphism part of the monoidal
structure. Therefore, consider $f:A\to A'$ and $g:B\to B'$ in
$\Nom/X$. We will define a morphism $f\otimes_X g: A\otimes_X B\to
A'\otimes_X B'$ as a family of maps $(f\otimes_X g)_x:(A\otimes_X B)_x\to
(A'\otimes_X B')_x$. They are defined by

\begin{center}
\begin{tikzpicture}
    \node (1) at (0,0) {$(A\otimes_X B)_x$};
    \node (2) at (4,0) {$(A'\otimes_X B')_x$};
    \node (3) at (0,-1.25) {$A_x\otimes_x B_x$}; 
    \node (4) at (4,-1.25) {$A'_x\otimes_x B'_x$.}; 
    
    \path[every node/.style={font=\sffamily\small}]
    (1) [-stealth,black,  thick] edge node [midway,above] {$(f\otimes_X g)_x$} (2)
    (3) [-stealth,black,  thick] edge node [midway,above] {$f_x\otimes_x g_x$} (4);
    
    \draw (1) [black,  thick] edge [double distance=1.5pt] (3);
    \draw (2) [black,  thick] edge [double distance=1.5pt] (4);
\end{tikzpicture}
\end{center}

To make sure that this is a well-defined morphism of nominal families,
by Remark~\ref{rk:NomFamPhiMap}, it is sufficient to show that these
maps are compatible with the family of maps $\eqfams{A\otimes_X B}{\pi}{x}$,
which we henceforth refer to as the $\varphi$-maps. In other words, we must
show that the outer part
of the following diagram commutes for all $\pi\in\PIA$, where
$y=\pi\bul x$.
\[
\begin{tikzpicture}[scale=0.8]
    \node (1) at (-8,0) {$\pihat(A_x\otimes_x B_x)$};             
    \node (2b) at (-3.5,-2) {$\pihat(A_{x})\otimes_{y} \pihat(B_{ x})$};       
    \node (3) at (-8,-4) {$A_{y}\otimes_{y} B_{y}$};
    \node (4) at (8,0) {$\pihat(A'_x\otimes_x B'_x)$};             
    \node (5b) at (3.5,-2) {$\pihat(A'_{x})\otimes_{y} \pihat(B'_{ x})$};       
    \node (6) at (8,-4) {$A'_{y}\otimes_{y} B'_{y}$.};

    \node (7) at (0,-0.75) {(1)};
    \node (8) at (-6.5,-2) {(2)};
    \node (9) at (6.5,-2) {(2)};
    \node (10) at (0,-2.75) {(3)};
    
    \path[every node/.style={font=\sffamily\small}]
	(2b) [-stealth,black,  thick] edge node [midway,right=0.8em]{$\eqfams{A}{\pi}{x} \otimes_{y} \eqfams{B}{\pi}{x}$} (3)
	(1) [-stealth,black,  thick] edge node [midway,left]{$\eqfams{A\otimes_X B}{\pi}{x}$} (3)
	
	(1) [-stealth,black,  thick] edge node [midway,above]{$\pihat(f_x\otimes_x g_x)$} (4)
	(1) [-stealth,black,  thick] edge node [midway,above]{$\mm_{\pi}$} (2b)
	(2b) [-stealth,black,  thick] edge node [midway,above]{$\pihat(f_x)\otimes_{y} \pihat(g_x)$} (5b)
	(3) [-stealth,black,  thick] edge node [midway,above]{$f_{y}\otimes_{y} g_{y}$} (6)
	(4) [-stealth,black,  thick] edge node [midway,above]{$\mm_{\pi}$} (5b)
	(5b) [-stealth,black,  thick] edge node [midway,left=1em]{$\eqfams{A'}{\pi}{x} \otimes_{y} \eqfams{B'}{\pi}{x}$} (6)
	(4) [-stealth,black,  thick] edge node [midway,right]{$\eqfams{A'\otimes_X B'}{\pi}{x}$} (6);
\end{tikzpicture}
\]
Part (1) commutes because $\mm_{\pi}$ is a natural transformation; the
parts labelled (2) commute by definition of the $\varphi$-maps for
$\otimes_X$; part (3) commutes since $f$ and $g$ are morphisms of
nominal families, and are therefore compatible with their respective
$\varphi$-maps.  This finishes the proof that $f\otimes_X g:A\otimes_X
B\to A'\otimes_X B'$ is a well-defined morphism.

Next, we must show that $\otimes_X$ is a functor. However, this is
straightforward because the required equations hold fiberwise: we have
$\id^A_x \otimes_x \id^B_x = \id^{A \otimes B}_x$ and $(f_x\otimes_x
g_x)\circ (f'_x\otimes_x g'_x)=(f_x\circ f'_x) \otimes_x (g_x\circ
g'_x)$, because $\otimes_x$ is a functor on $\NomAX{x}$ for every
$x\in X$.

\subsubsection{Definition of monoidal unit, associators, and unitors}
\label{sssec:def-unit-etc}

To continue the definition of the monoidal structure on $\Nom/X$, we
must provide a monoidal unit $I$, which we simply take to be the
nominal family $(I_x)_{x\in X}$, where $I_x$ is the monoidal unit of
$\otimes_x$ in $\NomAX{x}$. This becomes a nominal family via
$\varphi_{\pi,x} = \mm_{\pi}^0 : \pihat(I_x)\to I_{\pi\bul x}$.

We must also provide associators $\alpha^{A,B,C}:(A\otimes_X B)\otimes_X C\to
A\otimes_X(B\otimes_X C)$ and unitors $\lambda^{A}:I\otimes_X A\to A$ and $\rho^{A}:A\otimes_X
I\to A$.  We define these maps fiberwise, by taking each component
$\alpha_x$, $\lambda_x$, and $\rho_x$ to be an associator or unitor of
the corresponding monoidal structure $\otimes_x$ on $\NomAX{x}$:
\[
\begin{array}{lll}
  \alpha^{A,B,C}_{x} &:& (A_x\otimes_x B_x)\otimes_x C_x\to A_x\otimes_x(B_x\otimes_x C_x), \\
  \lambda^{A}_{x} &:& I_x\otimes_x A_x\to A_x,\\
  \rho^{A}_{x} &:& A_x\otimes_x I_x\to A_x.
\end{array}
\]
We must show that these maps are indeed morphisms of nominal
families. By Remark~\ref{rk:NomFamPhiMap}, it suffices to show that
they are compatible with the respective $\varphi$-maps.
This is shown in Appendix~\ref{app:associativity-etc}.
Then $\alpha^{A,B,C}$, $\lambda^{A}$, and
$\rho^{A}$ are well-defined morphisms of $\Nom/X$.

\subsubsection{Naturality and coherence}

To finish our definition of the monoidal structure on $\Nom/X$, we
must show that $\alpha$, $\lambda$, and $\rho$ are natural
transformations, and that they satisfy the required coherence diagrams
for monoidal categories, i.e., one pentagon and one triangle. However,
all of these properties are trivial because they hold separately in
each fiber.

This finishes the proof of Theorem~\ref{thm:monoidal-nom-slash-x}.

\section{$\Nom$ is locally monoidal closed}

Our goal is to show that the category $\Nom$ is locally monoidal
closed. More specifically, our main theorem states that every monoidal
closed structure on $\Nom$ induces a monoidal closed
structures on $\Nom/X$. 

\begin{xtheorem}\label{NomXFuncSpcObjThm}
  Every monoidal closed structure on $\Nom$ induces a
  monoidal closed structure on $\Nom/X$.
\end{xtheorem}

Before we prove this theorem, we first recall some basic facts about
monoidal functors between monoidal closed categories.

\subsection{Monoidal functors between monoidal closed categories}
\label{ssec:nmap-general}

Let $(\Cc,\otimes,\multimap)$ and $(\Dd,\otimes',\multimap')$ be
monoidal closed categories. It is well-known that a monoidal functor
$F:\Cc\to\Dd$ induces a natural transformation $\nn^{A,B}:F(A\multimap
B)\to FA \multimap' FB$. Specifically, $\nn^{A,B}$ is the unique map
making the following diagram commute:
\[
\begin{tikzpicture}
    \node (1) at (0,0) {$F(A\multimap B)\otimes' F(A)$};
    \node (2) at (0,-1.5) {$(F(A)\multimap' F(B))\otimes' F(A)$};
    \node (3) at (7,0) {$F((A\multimap B)\otimes A)$}; 
    \node (4) at (7,-1.5) {$F(B)$}; 
    
    \path[every node/.style={font=\sffamily\small}]
    (1) [-stealth,black,  thick] edge node [midway,left] {$\nn^{A,B}\otimes' \id$} (2)
    (1) [-stealth,black,  thick] edge node [midway,above] {$(\mm^{A\multimap B,A})\inv$} (3)
    (2) [-stealth,black,  thick] edge node [midway,below] {$\eps'$}  (4)
    (3) [-stealth,black,  thick] edge node [midway,right] {$F(\eps)$} (4);
\end{tikzpicture}
\]
Here, $\mm$ is the natural isomorphism making $F$ into a monoidal
functor, and $\eps$ and $\eps'$ are the application morphisms for the
respective monoidal closed structures. Equivalently, $\nn^{A,B} =
(F(\eps)\circ m_{A,B})^c$, where $(-)^c:(A\otimes'
B,C)\xrightarrow{\cong} (A,B\multimap' C)$ is the currying operation.
Moreover, if $F$ is isomorphism of categories, then $\nn^{A,B}$ is an
isomorphism.

\begin{xremark}\label{rmk:n-map-comp}
  If $(\Cc,\otimes,\multimap)$, $(\Dd,\otimes',\multimap')$, and
  $(\Ee,\otimes'',\multimap'')$ are monoidal closed categories and
  $F:\Cc\to\Dd$ and $G:\Dd\to\Ee$ are monoidal functors, the following
  diagram commutes.
  \[
  \begin{tikzpicture}
    \node (1) at (0,0) {$GF(B\multimap C)$};
    \node (2) at (3,-1.5) {$G(FB\multimap' FC)$};
    \node (3) at (6,0) {$GFB \multimap'' GFC$}; 
    
    \path[every node/.style={font=\sffamily\small}]
    (1) [-stealth,black,  thick] edge node [midway,left=0.5em] {$G(\nn_F)$} (2)
    (1) [-stealth,black,  thick] edge node [midway,above] {$\nn_{GF}$} (3)
    (2) [-stealth,black,  thick] edge node [midway,right=0.5em] {$\nn_G$}  (3);
  \end{tikzpicture}
  \]
  We also have $\nn_{\id} = \id$, where $\id:\Cc\to\Cc$ is the identity
  monoidal functor.
\end{xremark}

\subsection{Proof of Theorem~\ref{NomXFuncSpcObjThm}}

Fix an atom set $\aA$ and let $\Nom=\NomA$ be the category of nominal
sets over $\aA$. Let $X$ be an object of $\Nom$. Consider any monoidal
closed structure $(\otimes,\multimap)$ on $\Nom$. We wish to define a
monoidal closed structure $(\otimes_X,\multimap_X)$ on $\Nom/X$.
Since $\Nom/X$ is equivalent to the category $\Nom^X$ of nominal
families over $X$, it suffices to define a monoidal closed structure
on the latter.

The tensor $\otimes_X$ was already defined in
Section~\ref{sec:monoidal}; recall that its definition required some
fixed, but arbitrary choices of bijections $\delta_x:\aA\to\aA_x$, but
the resulting tensor was shown to be well-defined up to isomorphism.

We now define the closed structure in several steps.

\subsubsection{Definition of $A\multimap_X B$}\label{sssec:def-a-magicwand-b}

Consider any two objects $A,B\in\Nom^X$. We must define an object
$A\multimap_X B$ of $\Nom^X$. By definition, $A$ and $B$ are nominal
families, say $A = (A_x,\eqfams{A}{\pi}{x})$ and $B =
(B_x,\eqfams{B}{\pi}{x})$, where $A_x,B_x\in\NomAX{x}$. We wish to
define a nominal family $A\multimap_X B = (C_x)_{x\in X}$.

From Section~\ref{sec:monoidal}, we already have a monoidal structure
$\otimes_x$ on each fiber $\NomAX{x}$. Moreover, since this monoidal
structure was induced by some isomorphism with $\Nom$, it is closed.
Let $\multimap_x$ be the internal hom-functor on $\NomAX{x}$. Then we
can define the family of objects $C_x$ as follows:
\[
C_x = A_x\multimap_x B_x.
\]
To make this into a nominal family, we must define a family of
equivariant maps $\eqfams{C}{\pi}{x} : \pihat(C_x)\to C_{\pi\bul
  x}$, satisfying properties (\ref{ax-c}) and (\ref{ax-d}) of
Remark~\ref{rk:NomFamPhiMap}. We define these as follows. Here $y =
\pi\bul x$ as usual, and $\nn_{\pi}:\pihat(A_x\multimap_x
B_x)\to\pihat(A_x)\multimap_{y}\pihat(B_x)$ is the map induced from
$\mm_{\pi}:\pihat(A_x\otimes_x
B_x)\to\pihat(A_x)\otimes_{y}\pihat(B_x)$ as in
Section~\ref{ssec:nmap-general}.
\[
  \begin{tikzpicture}
    \node (1) at (0,0) {$\pihat(A_x\multimap_x B_x)$};
    \node (2) at (6,0) {$A_y\multimap_y B_y$}; 
    \node (3) at (3,-1.5) {$\pihat(A_x)\multimap_y\pihat(B_x)$};
    
    \path[every node/.style={font=\sffamily\small}]
    (1) [-stealth,black,  thick] edge node [midway,above] {$\eqfams{C}{\pi}{x}$} (2)
    (1) [-stealth,black,  thick] edge node [midway,left=0.5em] {$\nn_{\pi}$} (3)
    (3) [-stealth,black,  thick] edge node [midway,right=0.5em] {$(\eqfams{A}{\pi}{x})\inv \multimap_{y}\eqfams{B}{\pi}{x}$}  (2);
  \end{tikzpicture}
\]
Properties~(\ref{ax-c}) and (\ref{ax-d}) of
Remark~\ref{rk:NomFamPhiMap} are shown in Appendix~\ref{app:prop-a-b-2}.
This finishes the proof that $(C_x,\eqfams{C}{\pi}{x})$ is a nominal
family, i.e., that $A\multimap_X B$ is a well-defined object of
$\Nom^X$.

\subsubsection{Definition of $f\multimap_X g$ and functoriality}

To continue our definition of the closed structure on $\Nom/X$, we
must now define the morphism part of the functor $\multimap_X $.
Therefore, consider $f:A'\to A$ and $g:B\to B'$ in $\Nom/X$.  We will
define a morphism $f\multimap_X g:A\multimap_X B\to A'\multimap_X B'$
as a family of maps $(f\multimap_X g)x : (A\multimap_X
B)_x\to(A'\multimap_X B')_x$. They are defined by
\[
\begin{tikzpicture}
    \node (1) at (0,0) {$(A\multimap_X B)_x$};
    \node (2) at (4,0) {$(A'\multimap_X B')_x$};
    \node (3) at (0,-1.25) {$A_x\multimap_x B_x$}; 
    \node (4) at (4,-1.25) {$A'_x\multimap_x B'_x$.}; 
    
    \path[every node/.style={font=\sffamily\small}]
    (1) [-stealth,black,  thick] edge node [midway,above] {$(f\multimap_X g)_x$} (2)
    (3) [-stealth,black,  thick] edge node [midway,above] {$f_x\multimap_x g_x$} (4);
    
    \draw (1) [black,  thick] edge [double distance=1.5pt] (3);
    \draw (2) [black,  thick] edge [double distance=1.5pt] (4);
\end{tikzpicture}
\]
To make sure that this is a well-defined morphism of nominal families,
by Remark~\ref{rk:NomFamPhiMap}, it is sufficient to show that these
maps are compatible with the $\varphi$-maps, i.e., that the outer part
of the following diagram commutes for all $\pi\in\PIA$, where
$y=\pi\bul x$.
\[
\begin{tikzpicture}[scale=0.8]
    \node (1) at (-8,0) {$\pihat(A_x\multimap_x B_x)$};             
    \node (2b) at (-3.5,-2) {$\pihat(A_{x})\multimap_{y} \pihat(B_{ x})$};
    \node (3) at (-8,-4) {$A_{y}\multimap_{y} B_{y}$};
    \node (4) at (8,0) {$\pihat(A'_x\multimap_x B'_x)$};             
    \node (5b) at (3.5,-2) {$\pihat(A'_{x})\multimap_{y} \pihat(B'_{ x})$};
    \node (6) at (8,-4) {$A'_{y}\multimap_{y} B'_{y}$.};

    \node (7) at (0,-0.75) {(1)};
    \node (8) at (-6.5,-2) {(2)};
    \node (9) at (6.5,-2) {(2)};
    \node (10) at (0,-2.75) {(3)};
    
    \path[every node/.style={font=\sffamily\small}]
	(2b) [-stealth,black,  thick] edge node [midway,right=0.8em]{$(\eqfams{A}{\pi}{x})\inv \multimap_{y} \eqfams{B}{\pi}{x}$} (3)
	(1) [-stealth,black,  thick] edge node [midway,left]{$\eqfams{A\multimap_X B}{\pi}{x}$} (3)
	
	(1) [-stealth,black,  thick] edge node [midway,above]{$\pihat(f_x\multimap_x g_x)$} (4)
	(1) [-stealth,black,  thick] edge node [midway,above]{$\nn_{\pi}$} (2b)
	(2b) [-stealth,black,  thick] edge node [midway,above]{$\pihat(f_x)\multimap_{y} \pihat(g_x)$} (5b)
	(3) [-stealth,black,  thick] edge node [midway,above]{$f_{y}\multimap_{y} g_{y}$} (6)
	(4) [-stealth,black,  thick] edge node [midway,above]{$\nn_{\pi}$} (5b)
	(5b) [-stealth,black,  thick] edge node [midway,left=1em]{$(\eqfams{A'}{\pi}{x})\inv \multimap_{y} \eqfams{B'}{\pi}{x}$} (6)
	(4) [-stealth,black,  thick] edge node [midway,right]{$\eqfams{A'\multimap_X B'}{\pi}{x}$} (6);
\end{tikzpicture}
\]
Part (1) commutes because $\nn_{\pi}$ is a natural transformation; the
parts labelled (2) commute by definition of the $\varphi$-maps for
$\multimap_X$; part (3) commutes since $f$ and $g$ are morphisms of
nominal families, and are therefore compatible with their respective
$\varphi$-maps.  This finishes the proof that $f\multimap_X g:A\multimap_X
B\to A'\multimap_X B'$ is a well-defined morphism.

Next, we must show that $\multimap_X$ is a functor. However, this is
straightforward because the required equations hold fiberwise: we have
$\id^A_x \multimap_x \id^B_x = \id^{A \multimap B}_x$ and
$(f_x\multimap_x g_x)\circ (f'_x\multimap_x g'_x)=(f'_x\circ f_x)
\multimap_x (g_x\circ g'_x)$, because $\multimap_x$ is a (mixed
variance) functor on $\NomAX{x}$ for every $x\in X$.

\subsubsection{Adjunction between $\multimap_X$ and $\otimes_X$}

To finish our definition of the closed structure on $\Nom/X$, we must
show that the functor $B\multimap_X (-)$ is a right adjoint of
$(-)\otimes_X B$. Most of the proof is fiberwise, i.e., we already
know that the functor $B_x\multimap_x(-)$ is a right adjoint of
$(-)\otimes_x B_x$ in every fiber $x\in X$.  The only additional thing
we must prove is that the fiberwise adjunction respects morphisms of
nominal families, i.e., that if $f_x:A_x\otimes_x B_x\to C_x$ and
$g_x:A_x\to B_x\multimap_x C_x$ are fiberwise adjoint mates, then the
family $(f_x)_{x\in X}$ respects the $\varphi$-maps if and only if the
family $(g_x)_{x\in X}$ does.

So consider a family of maps $f_x:A_x\otimes_x B_x\to C_x$, and let
$g_x:A_x\to B_x\multimap_x C_x$ be the curry of $f_x$ for every $x\in
X$. We must show that the perimeter of the left diagram commutes if
and only if the perimeter of the right one does:
\[
\begin{tikzpicture}
    
    \node (lltop) at (0,0) {$\pihat(A_x\otimes_x B_x)$};
    \node (lrtop) at (3.5,0) {$\pihat(C_x)$};
    \node (llbot) at (0,-4) {$A_{\pi \bul x}\otimes_{\pi \bul x} B_{\pi \bul x}$}; 
    \node (lrbot) at (3.5,-4) {$C_{\pi \bul x}$}; 
    \node (llmid) at (0,-2) {$\pihat(A_x)\otimes_{\pi \bul x} \pihat(B_x)$};
    
    \node (rltop) at (6,0) {$\pihat(A_x)$};
    \node (rrtop) at (9.5,0) {$\pihat(B_x\multimap_x C_x)$};
    \node (rlbot) at (6,-4) {$A_{\pi \bul x}$}; 
    \node (rrbot) at (9.5,-4) {$B_{\pi \bul x} \multimap_{\pi \bul x} C_{\pi \bul x}$}; 
    \node (rrmid) at (9.5,-2) {$\pihat(B_x)\multimap_{\pi \bul x} \pihat(C_x)$};

    \node (rlabel0) at (-0.75,-1) {(0)};
    \node (rlabel1) at (7,-2) {(3)};
    \node (rlabel2) at (8.5,-0.75) {(4)};
    \node (llabel0) at (10.25,-1) {(5)};
    \node (llabel1) at (2.5,-2) {(2)};
    \node (llabel2) at (1.3,-0.75) {(1)};

    \path[every node/.style={font=\sffamily\small}]
	(rltop) [-stealth,black,  thick] edge node [midway,above]{$\pihat(g_x)$} (rrtop) 
	(rltop) [-stealth,black,  thick] edge node [midway,left]{$\eqfams{A}{\pi}{x}$} (rlbot) 
	(rlbot) [-stealth,black,  thick] edge node [midway,below]{$g_{\pi\bul x}$} (rrbot)
	(rrtop) [-stealth,black,  thick] edge node [midway,left]{$\nn_{\pi}$} (rrmid)
	(rrmid) [-stealth,black,  thick] edge node [midway,left]{$\left(\eqfams{B}{\pi}{x}\right)\inv\multimap_{\pi \bul x} \eqfams{C}{\pi}{x}$} (rrbot)
	(rltop) [-stealth,black,  thick] edge node [midway,left]{$h^c\:\:$} (rrmid)
	
	(rrtop) [-stealth,black,  thick] edge [out= -20, in=20] node [midway,right]{$\eqfams{B\multimap_X C}{\pi}{x}$} (rrbot)    	    
    
	(lltop) [-stealth,black,  thick] edge node [midway,above]{$\pihat(f_x)$} (lrtop) 
	(llbot) [-stealth,black,  thick] edge node [midway,below]{$f_{\pi\bul x}$} (lrbot)
	(lrtop) [-stealth,black,  thick] edge node [midway,right]{$\eqfams{C}{\pi}{x}$} (lrbot)
	(lltop) [-stealth,black,  thick] edge node [midway,right]{$\mm_{\pi}$} (llmid)
	(llmid) [-stealth,black,  thick] edge node [midway,right]{$\eqfams{A}{\pi}{x}\otimes_{\pi \bul x} \eqfams{B}{\pi}{x}$} (llbot)
	(llmid) [-stealth,black,  thick] edge node [midway,right]{$\:\:h$} (lrtop)
	
	(lltop) [-stealth,black,  thick] edge [out = 200, in= 160] node [midway,left]{$\eqfams{A\otimes_X B}{\pi}{x}$} (llbot) ;
\end{tikzpicture}
\]
To see why, first define $h$ to be the unique map such that (1)
commutes. Let $h^c$ be the result of currying $h$; then (4) commutes
by the definition of $\nn$ and the fact that $g_x$ if the curried form
of $f_x$. Parts (0) and (5) commute by the definitions of
$\eqfams{A\otimes_X B}{\pi}{x}$ and $\eqfams{B\multimap_X C}{\pi}{x}$,
respectively. Then by the naturality of currying, (2) commutes if and
only if (3) commutes, which implies that the perimeter of the left
diagram commutes if and only if the perimeter of the right one does,
as claimed.

This finishes the proof of Theorem~\ref{NomXFuncSpcObjThm}.

\section{Conclusion and future work}

We showed that every monoidal (closed) structure on $\Nom$ induces a
corresponding monoidal (closed) structure on all of its slice
categories $\Nom/X$. It seems reasonable to say that $\Nom$ is
therefore ``locally monoidal closed''. However, as far as we know,
there is no general definition of local monoidal structure: given a
general category $\Cc$, it may not in general make sense to ask
whether a monoidal structure on $\Cc/X$ is ``induced'' by a monoidal
structure on $\Cc$. It is an interesting question what requirements a
category $\Cc$ should satisfy so that one can speak of locally
monoidal and locally monoidal closed structures.

A possible limitation of our result is that there are only a few
interesting examples of monoidal closed structure on $\Nom$. The
best-known ones are the cartesian structure and the separated one;
another (non-symmetric) monoidal closed structure was very recently
discovered by \cite{lenke2026}. For each of these cases, the
corresponding monoidal closed structures on the slice categories could
have been defined in elementary terms, yielding a more concrete
element-wise proof. Nevertheless, we believe there is some value in
having proved this result in full generality. For one, it is pleasing
that the proof is uniform, i.e., it works identically regardless of
the particular monoidal closed structure in question. Second, part of
our result is not just for monoidal closed structures, but more
generally for monoidal ones, of which there are many. Third, the
technique used in this paper is clearly not limited to just monoidal
structures, but can be expected to work without much difficulty for
other similar categorical structures, i.e., those that are given by
functors, natural transformations, and equations. Thus, we expect that
there is a much larger class of structures that can be lifted from
$\Nom$ to its slice categories. The details of this are left for
future work.

As mentioned in the introduction, one of the motivations for studying
the structure of slice categories on $\Nom$ is that slice categories
play an important role in models of dependent type theory. We hope
that the present work can serve as a building block towards the
eventual goal of constructing a model for a nominal dependent type
theory. 


\bibliographystyle{./entics}
\bibliography{localnom}

\appendix

\section{Some proofs}\label{app:proofs}

\begin{proof*}{Proof of Proposition~\ref{prop:pi-extension}.}
Let $\sigma\in \PIL$ and $x\in X$. We would like to define $\sigma \bulL x
= \pi \bul x$, where $\pi$ is some finitely supported permutation
such that $\pi\restr{\supp(x)}=\sigma\restr{\supp(x)}$. For this to
work, we must show that:
\begin{enumerate}

\item Such a permutation $\pi$ exists. Indeed, we know $\supp(x)$ is
  finite. Let $\supp(x)=\set{a_1,a_2,\ldots, a_n}$. Let $b_i=\sigma(a_i)$
  for $i=1,\ldots,n$. Let $U=\set{a_1,a_2,\ldots, a_n} \cup
  \set{b_1,b_2,\ldots, b_n}$, which may have any number of elements
  between $n$ and $2n$, say $n+k$ elements. Let $\set{c_1,c_2,\ldots,
    c_k} = U\setminus\set{a_1,a_2,\ldots, a_n}$ and $\set{d_1,d_2,\ldots,
    d_k}=U\setminus\set{b_1,b_2,\ldots, b_n}$, and define $\pi\in \PI$ by
  \[
  \pi(a_i) = b_i \quad \mbox{and} \quad \pi(c_j) = d_j
  \]
  for $i=1,\ldots,n$ and $j=1,\ldots,k$. Also, by definition
  $\pi\restr{\supp(x)}=\sigma\restr{\supp(x)}$.

\item The choice of $\pi$ is not important, i.e., if $\pi$ and
  $\pi'$ satisfy
  $\pi\restr{\supp(x)}=\pi'\restr{\supp(x)}=\sigma\restr{\supp(x)}$,
  then $\pi \bul x = \pi' \bul x$. This holds by Lemma
  \ref{LemSuppAct}. Hence, $\bulL:\PIL\times X \to X$ is a
  well-defined operation. We also note that if $\sigma$ is finitely
  supported, we have $\sigma\bulL x = \sigma\bul x$.

\item The operation $\bulL$ is an action, i.e., $\id\bulL x =
  x$ and for all $\sigma_1,\sigma_2\in \PIL$, $(\sigma_2\sigma_1)\bulL x =
  \sigma_2\bulL (\sigma_1\bulL x)$. For the first claim, we clearly
  have $\id\bulL x = \id\bul x = x$. For the second claim,
  choose some $\pi_1,\pi_2\in\PI$ such that
  \begin{equation}\label{eqn:action-extend1}
    \pi_1\restr{\supp(x)}=\sigma_1\restr{\supp(x)}
  \end{equation}
  and
  \begin{equation}\label{eqn:action-extend2}
    \pi_2\restr{\supp(\pi_1\bul x)}=\sigma_2\restr{\supp(\pi_1\bul x)}.
  \end{equation}

  We first claim that
  $\pi_2\pi_1\restr{\supp(x)}=\sigma_2\sigma_1\restr{\supp(x)}$. Indeed,
  for $a\in \supp(x)$, we have $\sigma_2\sigma_1 (a) = \sigma_2 (\sigma_1 (a)) =
  \sigma_2 (\pi_1 (a)) = \pi_2 (\pi_1 (a)) = \pi_2 \pi_1 (a)$, where
  the second and third equations hold by {\eqref{eqn:action-extend1}}
  and {\eqref{eqn:action-extend2}}, respectively.  It follows that
  $(\sigma_1\sigma_2)\bulL x = (\pi_1\pi_2)\bul x$.  Then we have
  \[
  \begin{array}{rcl@{\qquad}l}
    \sigma_2\bulL (\sigma_1\bulL x)
    & = & \sigma_2\bulL (\pi_1\bul x)\\
    & = & \pi_2\bul (\pi_1\bul x)\\
    & = & (\pi_2\pi_1)\bul x \\
    & = & (\sigma_2\sigma_1)\bulL x,
  \end{array}
  \]
  as claimed. Thus, the operation $\bulL$ is a valid action.\qedhere
\end{enumerate}
\end{proof*}

\begin{proof*}{Proof of Proposition~\ref{DotBlankMapNatIsoProp}.}
  We first need to show that $\DBL{\sigma}:\sigmahat(X)\to X$ is a
  well-defined morphism, i.e., that it is equivariant. So consider any
  $\pi\in\PIA$ and $x\in X$. We have
  \[
  \sigma\bulL(\pi\diamond x)
  = \sigma\bulL((\sigma\inv\pi\sigma)\bul x)
  = \sigma\bulL\sigma\inv\bulL\pi\bulL\sigma\bulL x
  = \pi\bul(\sigma\bulL x).
  \]
  Next, we need to show naturality. So consider objects $(X,\bul)$,
  $(Y,\diamond)\in\NomA$.  Here, for clarity, we have used two
  different symbols for the actions on $X$ and $Y$. Let $f:X\to Y$ be
  an equivariant map.  We must show that the following square
  commutes:
  \[
  \begin{tikzpicture}
    \node (ltop) at (0,0) {$\sigmahat(X)$};
    \node (rtop) at (3,0) {$ \sigmahat(Y)$};
    \node (lbot) at (0,-1.5) {$X$}; 
    \node (rbot) at (3,-1.5) {$Y$}; 
    
    \draw [-stealth][black,thick] (ltop) -- (rtop) node[midway,above] {$\sigmahat(f)$};
    \draw [-stealth][black,thick] (ltop) -- (lbot) node[midway,left] {$\DBL{\sigma}$};
    \draw [-stealth][black,thick] (lbot) -- (rbot) node[midway,below] {$f$};
    \draw [-stealth][black,thick] (rtop) -- (rbot) node[midway,right] {$\DDL{\sigma}$};
  \end{tikzpicture}
  \]
  In other words, for all $x\in X$, we need to show that $f(\sigma\bulL
  x)=\sigma\diamondL f(x)$.  Let $\pi$ be a finitely supported permutation
  such that $\pi$ and $\sigma$ agree on the $\supp(x)$. Note that
  $\supp(f(x))\subseteq \supp(x)$, so $\pi$ and $\sigma$ agree on
  $\supp(f(x))$ as well. Then, because $f$ is equivariant, we have
  \[
  f(\sigma \bulL x)
  = f(\pi\bul x)
  = \pi \diamond f(x)
  = \sigma \diamondL f(x).
  \]
  Finally, the map $\DBL{\sigma}$ is clearly invertible with inverse
  $\DBL{\sigma\inv}$.   
\end{proof*}

\begin{proof*}{Proof of Lemma~\ref{lem:EtaCompose}.}
  Define $\sigma_1,\sigma_2,\sigma_3:\bB \to \bB$ by $\sigma_1=\gamma_2\circ
  (\gamma_1)\inv$, $\sigma_2=\gamma_3\circ (\gamma_2)\inv$, and
  $\sigma_3=\gamma_3\circ (\gamma_1)\inv$. Note that $\sigma_3 = \sigma_2\sigma_1$.
  Consider $\sigmahat_1,\sigmahat_2,\sigmahat_3:\NomB\to\NomB$ as in
  Section~\ref{ssec:muhat-bijective}. Consider the following diagram:
  \[
  \resizebox{\textwidth}{!}{
  \begin{tikzpicture}
    \node (1) at (-6,-3) {$A\otimes_1 B$}; 
    \node (2) at (0,-3) {$A\otimes_2 B$}; 
    \node (3) at (6,-3) {$A\otimes_3 B$}; 
    \node (4) at (-6,0) {$\widehat{\sigma_1}(A\otimes_1 B)$};
    \node (5) at (-3,0) {$\widehat{\sigma_1}(A)\otimes_2 \widehat{\sigma_1}(B)$};
    \node (6) at (3,0) {$\widehat{\sigma_2}(A\otimes_2 B)$};
    \node (7) at (6,0) {$\widehat{\sigma_2}(A)\otimes_3 \widehat{\sigma_2}(B)$};
    \node (8) at (-6,3) {$\widehat{\sigma_2}(\widehat{\sigma_1}(A\otimes_1 B))$};
    \node (9) at (0,3) {$\widehat{\sigma_2}(\widehat{\sigma_1}(A)\otimes_2 \widehat{\sigma_1}(B))$};
    \node (10) at (6,3) {$\widehat{\sigma_2}(\widehat{\sigma_1}(A))\otimes_3 \widehat{\sigma_2}(\widehat{\sigma_1}(B))$};
    \node (11) at (-3,4) {$\widehat{\sigma_3}(A\otimes_1 B)$};
    \node (12) at (3,4) {$\widehat{\sigma_3}(A)\otimes_3 \widehat{\sigma_3}(B)$};
    
    \node at (-7,1) {{\Large $(3)$}};  
    \node at (-3.5,1.5) {{\Large $(5)$}};  
    \node at (-3.5,-1.5) {{\Large $(1)$}};  
    \node at (0,0) {{\Large $(7)$}}; 
    \node at (3.5,1.5) {{\Large $(6)$}};  
    \node at (3.5,-1.5) {{\Large $(2)$}};  
    \node at (7,1) {{\Large $(4)$}};  
    \node at (0,-3.7) {{\Large $(8)$}}; 
    
    \path[every node/.style={font=\sffamily\small}]
    (4) [black,thick] edge [double distance=2pt] (5)
    (6) [black,thick] edge [double distance=2pt] (7)
    (8) [black,thick] edge [double distance=2pt] (9)
    (9) [black,thick] edge [double distance=2pt] (10)
    (8) [black,thick] edge [double distance=2pt] (11)
    (11) [black,thick] edge [double distance=2pt, bend left=8] (12)
    (12) [black,thick] edge [double distance=2pt] (10)
    
    (1) [-stealth,black,thick] edge node [midway,below]{$\eta_{\gamma_1,\gamma_2}^{A,B}$} (2)
    (2) [-stealth,black,thick] edge node [midway,below]{$\eta_{\gamma_2,\gamma_3}^{A,B}$} (3)
    (1) [-stealth,black,thick] edge [bend right=20] node [midway,below]{$\eta_{\gamma_1,\gamma_3}^{A,B}$} (3)
    
    (4) [-stealth,black,thick] edge node[midway,right] {$\DBL{\sigma_1}$} node[midway,left] {$\cong$} (1)
    (5) [-stealth,black,thick] edge node[near start,left] {$\cong$} node[near start,right] {$\:\:\DBL{\sigma_1}\otimes_2 \DBL{\sigma_1}$} (2)
    (6) [-stealth,black,thick] edge node[near end,right] {$\DBL{\sigma_2}$} node[near end,left] {$\cong$} (2)
    (7) [-stealth,black,thick] edge node[near start,right] {$\cong$} node[near start,left] {$\DBL{\sigma_2}\otimes_3 \DBL{\sigma_2}$} (3)
    
    (8) [-stealth,black,thick] edge node[midway,right] {$\DDL{\sigma_2}$} node[midway,left] {$\cong$} (4)
    (9) [-stealth,black,thick] edge node[near start,left] {$\DDL{\sigma_2}$} node[near start,right] {$\cong$} (5)
    (9) [-stealth,black,thick] edge node[near end,right] {$\cong$} node[near end,left] {$\widehat{\sigma_2} (\DBL{\sigma_1}\otimes_2 \DBL{\sigma_1})\:\:$} (6)
    (10) [-stealth,black,thick] edge node[near start,right] {$\cong$} node[near start,left] {$\widehat{\sigma_2} (\DBL{\sigma_1}) \otimes_3 \widehat{\sigma_2} (\DBL{\sigma_1})$} (7)
    
    (8) [-stealth,black,thick] edge [out=210, in=150, looseness=1.1] node[midway,left] {$\DBL{\sigma_3}$} node[midway,right] {$\cong$} (1)
    (10) [-stealth,black,thick] edge [out=-30, in=30, looseness=1.1] node[midway,right] {$\DBL{\sigma_3} \otimes_3 \DBL{\sigma_3}$} node[midway,left] {$\cong$} (3);
  \end{tikzpicture}
  }
  \]
  We must show that (8) holds. Note that (1), (2), and the perimeter
  commute by Remark~\ref{rem:EtaMap}. (3) and (4) commute by
  Lemma~\ref{lem:mu1mu2}. (5) and (6) commute because $\sigmahat_1$ and
  $\sigmahat_2$ strictly preserve the relevant monoidal structures, and
  (7) commutes by naturality of $\DBL{\sigma_2}$. Since all maps in this
  diagram are isomorphisms, it follows that (8) commutes, as claimed.
\end{proof*}  

\begin{proof*}{Proof of Lemma~\ref{lem:eta-pihat}.}
  Let $\sigma=\gamma_2\gamma_1\inv$. Note that $(\pi\gamma_2)(\pi\gamma_1)\inv =
  \pi\sigma\pi\inv$. For clarity, we use ``$\bul$'' to denote all of
  the actions in the category $\NomA$ and ``$\diamond$'' for the
  actions in $\NomB$. Consider the following diagram.
  \[
  \resizebox{\textwidth}{!}{
  \begin{tikzpicture}
    \node (1) at (0,2) {$\widehat{\pi}(\widehat{\sigma}(A\otimes_1 B))$};
    \node (2) at (5,2) {$\widehat{\pi}(\widehat{\sigma}(A)\otimes_2 \widehat{\sigma}(B))$};
    
    \node (3) at (-7,0) {$\widehat{\pi\sigma}(A\otimes_1 B)$};
    \node (4) at (0,0) {$\widehat{\pi}(A\otimes_1 B)$};
    \node (5) at (5,0) {$\widehat{\pi}(A\otimes_2 B)$};
    \node (6) at (12,0) {$\widehat{\pi}(\widehat{\sigma}(A)) \otimes'_2 \widehat{\pi}(\widehat{\sigma}(B))$};
    
    \node (7) at (-7,-2) {$\widehat{\pi\sigma\pi\inv}(\widehat{\pi}(A\otimes_1 B))$};
    \node (8) at (0,-2) {$\widehat{\pi}(A)\otimes'_1 \widehat{\pi}(B)$};
    \node (9) at (5,-2) {$\widehat{\pi}(A)\otimes'_2 \widehat{\pi}(B)$};
    \node (10) at (12,-2) {$\widehat{\pi\sigma}(A)\otimes'_2 \widehat{\pi\sigma}(B)$};
    
    \node (11) at (-5,-4) {$\widehat{\pi\sigma\pi\inv}(\widehat{\pi}(A)\otimes'_1 \widehat{\pi}(B))$};
    \node (12) at (10,-4) {$\widehat{\pi\sigma\pi\inv}(\widehat{\pi}(A))\otimes'_2 \widehat{\pi\sigma\pi\inv}(\widehat{\pi}(B))$};

    \node at (2.5,-1) {$(1)$};
    \node at (1,1) {$(2)$};
    \node at (1,-3) {$(2)$};
    \node at (-3,-2) {$(3)$};
    \node at (8,0.5) {$(4)$};
    \node at (9,-2) {$(5)$};
    \node at (-4,0.5) {$(5)$};

    \path[every node/.style={font=\sffamily\small}]
    (1) [black,thick] edge [double distance=2pt] (2)
    (1) [black,thick] edge [double distance=2pt,bend right=20] (3)
    (2) [black,thick] edge [double distance=2pt,bend left=20] (6)
    (3) [black,thick] edge [double distance=2pt] (7)
    (4) [black,thick] edge [double distance=2pt] (8)
    (5) [black,thick] edge [double distance=2pt] (9)
    (6) [black,thick] edge [double distance=2pt] (10)
    (7) [black,thick] edge [double distance=2pt,bend right=20] (11)
    (10) [black,thick] edge [double distance=2pt,bend left=20] (12)
    (11) [black,thick] edge [double distance=2pt] (12)
    
    (1) [-stealth,black,thick] edge node [midway,left]{$\widehat{\pi}(\DBL{\sigma})$} (4)
    (2) [-stealth,black,thick] edge node [midway,left]{$\widehat{\pi}(\DBL{\sigma} \otimes_2 \DBL{\sigma})$} (5)
    (4) [-stealth,black,thick] edge node [midway,below]{$\widehat{\pi}(\eta_{\gamma_1,\gamma_2}^{A,B})$} (5)
    (6) [-stealth,black,thick] edge node [near start,left]{$\widehat{\pi}(\DBL{\sigma}) \otimes'_2 \widehat{\pi}(\DBL{\sigma})\:\:$} (9)
    (8) [-stealth,black,thick] edge node [midway,below]{$\eta_{\pi\gamma_1,\pi\gamma_2}^{A,B}$} (9)
    (7) [-stealth,black,thick] edge node [midway,above]{$\DDL{\pi\sigma\pi\inv}\:\:$} (4)
    (11) [-stealth,black,thick] edge node [midway,right]{$\:\:\DDL{\pi\sigma\pi\inv}$} (8)
    (12) [-stealth,black,thick] edge node [midway,left]{$\DDL{\pi\sigma\pi\inv} \otimes'_2 \DDL{\pi\sigma\pi\inv}$} (9);
  \end{tikzpicture}
  }
  \]
  We must show that (1) commutes. The parts labelled (2) commute by
  definition of $\eta^{A,B}$ map, (3) commutes by equality of nominal
  sets, (4) commutes since $\widehat{\pi}$ is a monoidal functor, and
  the two parts labelled (5) commute by
  Lemma~\ref{lem:permutationconjugatemap}.
\end{proof*}

\begin{proof*}{Proof of Lemma~\ref{lem:mmapcompose}.}
  On $\NomC$, let $\otimes'_2$ be the monoidal structure induced by
  $\mu_1\delta_1$. On $\NomD$, let $\otimes'_3$ and $\otimes''_3$ be
  the monoidal structures induced by $\mu_2\delta_2$ and
  $\mu_2\mu_1\delta_1$, respectively. Consider the following diagram. 
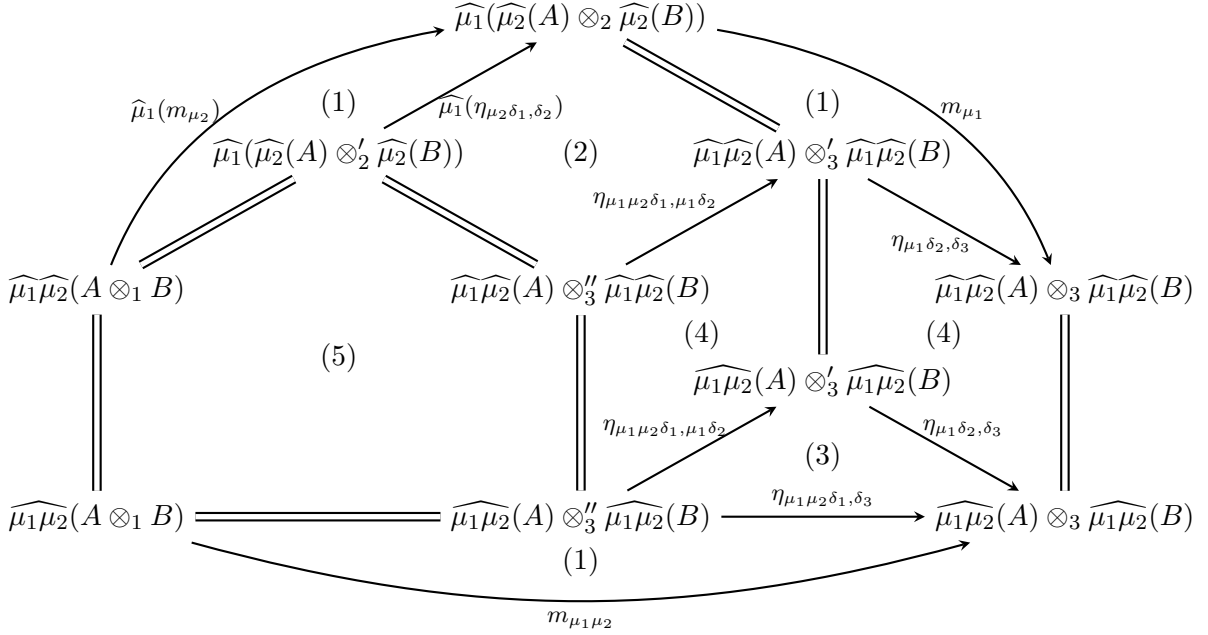
\[
\begin{tikzpicture}[xscale=0.8,yscale=0.6]
    \node (1) at (0,5) {$\widehat{\mu_1}\widehat{\mu_2}(A\otimes_1 B)$};     
    \node (2) at (8,11) {$\widehat{\mu_1}(\widehat{\mu_2}(A)\otimes_2 \widehat{\mu_2}(B))$};      
    \node (1a) at (4,8) {$\widehat{\mu_1}(\widehat{\mu_2}(A)\otimes'_2 \widehat{\mu_2}(B))$};      
    
    \node (4) at (0,0) {$\widehat{\mu_1\mu_2}(A\otimes_1 B)$};
    \node (5) at (8,0) {$\widehat{\mu_1\mu_2}(A)\otimes''_3 \widehat{\mu_1\mu_2}(B)$};
    \node (6) at (12,8) {$\widehat{\mu_1}\widehat{\mu_2}(A)\otimes'_3 \widehat{\mu_1}\widehat{\mu_2}(B)$};
    \node (5a) at (12,3) {$\widehat{\mu_1\mu_2}(A)\otimes'_3 \widehat{\mu_1\mu_2}(B)$};     
    \node (6a) at (8,5) {$\widehat{\mu_1}\widehat{\mu_2}(A)\otimes''_3 \widehat{\mu_1}\widehat{\mu_2}(B)$};
    
    \node (8) at (16,0) {$\widehat{\mu_1\mu_2}(A)\otimes_3 \widehat{\mu_1\mu_2}(B)$};
    \node (9) at (16,5) {$\widehat{\mu_1}\widehat{\mu_2}(A)\otimes_3 \widehat{\mu_1}\widehat{\mu_2}(B)$};
    
    \node at (8,8) {$(2)$};
    \node at (10,4) {$(4)$};
    \node at (12,1.33) {$(3)$};
    \node at (14,4) {$(4)$};
    \node at (4,9.1) {$(1)$};
    \node at (12,9.1) {$(1)$};
    \node at (8,-1) {$(1)$};
    \node at (4,3.5) {$(5)$};
    
    \path[every node/.style={font=\sffamily\small}]
    (1) [black,thick] edge [double distance=2pt] (4)
    (1) [black,thick] edge [double distance=2pt] (1a)
    (4) [black,thick] edge [double distance=2pt] (5)
    (1a) [black,thick] edge [double distance=2pt] (6a)
    (5) [black,thick] edge [double distance=2pt] (6a)
    (5a) [black,thick] edge [double distance=2pt] (6)
    (2) [black,thick] edge [double distance=2pt] (6)
    (8) [black,thick] edge [double distance=2pt] (9)
    
    (1a) [-stealth,black,  thick] edge node [near start,right] {$\:\widehat{\mu_1} (\eta_{\mu_2 \delta_1, \delta_2})\:$} (2)
    (5) [-stealth,black,  thick] edge node [midway,above] {$\eta_{\mu_1\mu_2 \delta_1, \delta_3}$} (8)
    (5) [-stealth,black,  thick] edge node [near end,left] {$\eta_{\mu_1\mu_2 \delta_1, \mu_1\delta_2}$} (5a)
    (6a) [-stealth,black,  thick] edge node [near end,left] {$\eta_{\mu_1\mu_2 \delta_1, \mu_1\delta_2}\:\:$} (6)
    (5a) [-stealth,black,  thick] edge node [near start,right] {$\:\eta_{\mu_1 \delta_2, \delta_3}$} (8)
    
    (6) [-stealth,black,  thick] edge node [near end,left] {$\eta_{\mu_1\delta_2, \delta_3}$} (9)  
    (1) [-stealth,black,  thick] edge [bend left=30] node [midway,left] {$\muhat_1(\mm_{\mu_2})$} (2)
    (2) [-stealth,black,  thick] edge [bend left=30] node [midway,right] {$\mm_{\mu_1}$} (9)
    (4) [-stealth,black,  thick] edge [bend right=20] node [midway,below] {$\mm_{\mu_1\mu_2}$} (8)
    ;
\end{tikzpicture}
\]
We must show that the outside commutes. Each part labelled (1)
commutes by definition of the $\mm$-maps; (2) commutes by
Lemma~\ref{lem:eta-pihat}; (3) commutes by Lemma~\ref{lem:EtaCompose};
each part labelled (4) commutes by naturality of $\eta$; and (5)
commutes because all its bounding morphisms are equalities.
\end{proof*}

\begin{proof*}{Proof of Lemma~\ref{lem:m-pi}.}
  Consider the following diagram, where $\sigma =
  \delta_2(\mu\delta_1)\inv = \mu\inv$.
  \[
  \begin{tikzpicture}
    \node (1) at (0,0) {$\muhat(A\otimes_1 B)$};
    \node (2) at (5,0) {$\muhat(A)\otimes_2'\muhat(B)$};
    \node (3) at (10,0) {$\muhat(A)\otimes_2\muhat(B)$};
    \node (4) at (0,-1.5) {$\sigmahat\muhat(A\otimes_1 B)$};
    \node (5) at (5,-1.5) {$\sigmahat(\muhat(A)\otimes_2'\muhat(B))$};
    \node (6) at (10,-1.5) {$\sigmahat\muhat(A)\otimes_2\sigmahat\muhat(B)$};
    \node (7) at (0,-2.75) {$A\otimes_1 B$};
    \node (9) at (10,-2.75) {$A\otimes_2 B$.};

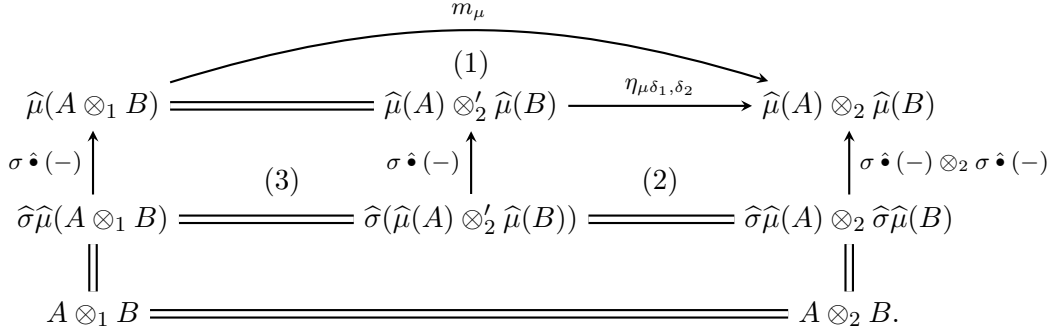
    \path[every node/.style={font=\sffamily\small}]
    (4) [black,thick] edge [double distance=2pt] (5)
    (5) [black,thick] edge [double distance=2pt] (6)
    (1) [black,thick] edge [double distance=2pt] (2)
    (4) [black,thick] edge [double distance=2pt] (7)
    (6) [black,thick] edge [double distance=2pt] (9)
    (7) [black,thick] edge [double distance=2pt] (9)
    (1) [-stealth,black,thick] edge [bend left=17] node [midway,above] {$\mm_{\mu}$} (3)
    (2) edge node [midway,above] {$\eta_{\mu\delta_1,\delta_2}$} (3)
    (5) [black,thick] edge node [midway,left]{$\DBL{\sigma}$} (2)
    (4) [black,thick] edge node [midway,left]{$\DBL{\sigma}$} (1)
    (6) [black,thick] edge node [midway,right] {$\DBL{\sigma}\otimes_2\DBL{\sigma}$} (3)
    ;

    \node at (5,0.55) {$(1)$};
    \node at (7.5,-1) {$(2)$};
    \node at (2.5,-1) {$(3)$};
  \end{tikzpicture}
  \]
  The part labelled (1) commutes by definition of $\mm_{\mu}$. The part
  labelled (2) commutes by Remark~\ref{rem:EtaMap}. The part labelled (3)
  commutes because the left and right objects are literally the same
  nominal set. Therefore, the outside commutes. This implies the Lemma
  because $\mu=\sigma\inv$. 
\end{proof*}

\section{Proof of properties (a) and (b) from Section~\ref{sssec:def-a-times-b}}
\label{app:prop-a-and-b-1}

\begin{xlem}\label{lem:prop-a-and-b-1}
  The family of maps $\eqfams{C}{\pi}{x}$ defined in
  Section~\ref{sssec:def-a-times-b} satisfies properties (\ref{ax-c})
  and (\ref{ax-d}) of Remark~\ref{rk:NomFamPhiMap}.
\end{xlem}

\begin{proof}
  We first need to show that $\DBL{\sigma}:\sigmahat(X)\to X$ is a
  well-defined morphism, i.e., that it is equivariant. So consider any
  $\pi\in\PIA$ and $x\in X$. We have
  \[
  \sigma\bulL(\pi\diamond x)
  = \sigma\bulL((\sigma\inv\pi\sigma)\bul x)
  = \sigma\bulL\sigma\inv\bulL\pi\bulL\sigma\bulL x
  = \pi\bul(\sigma\bulL x).
  \]
  Next, we need to show naturality. So consider objects $(X,\bul)$,
  $(Y,\diamond)\in\NomA$.  Here, for clarity, we have used two
  different symbols for the actions on $X$ and $Y$. Let $f:X\to Y$ be
  an equivariant map.  We must show that the following square
  commutes:
  \[
  \begin{tikzpicture}
    \node (ltop) at (0,0) {$\sigmahat(X)$};
    \node (rtop) at (3,0) {$ \sigmahat(Y)$};
    \node (lbot) at (0,-1.5) {$X$}; 
    \node (rbot) at (3,-1.5) {$Y$}; 
    
    \draw [-stealth][black,thick] (ltop) -- (rtop) node[midway,above] {$\sigmahat(f)$};
    \draw [-stealth][black,thick] (ltop) -- (lbot) node[midway,left] {$\DBL{\sigma}$};
    \draw [-stealth][black,thick] (lbot) -- (rbot) node[midway,below] {$f$};
    \draw [-stealth][black,thick] (rtop) -- (rbot) node[midway,right] {$\DDL{\sigma}$};
  \end{tikzpicture}
  \]
  In other words, for all $x\in X$, we need to show that $f(\sigma\bulL
  x)=\sigma\diamondL f(x)$.  Let $\pi$ be a finitely supported permutation
  such that $\pi$ and $\sigma$ agree on the $\supp(x)$. Note that
  $\supp(f(x))\subseteq \supp(x)$, so $\pi$ and $\sigma$ agree on
  $\supp(f(x))$ as well. Then we have
  \[
  f(\sigma \bulL x)
  = f(\pi\bul x)
  = \pi \diamond f(x)
  = \sigma \diamondL f(x).
  \]
  Finally, the map $\DBL{\sigma}$ is clearly invertible with inverse
  $\DBL{\sigma\inv}$.   
\end{proof}

\section{Proof of compatibility of the associators and unitors from Section~\ref{sssec:def-unit-etc}}
\label{app:associativity-etc}

Recall that we defined
\[
\begin{array}{lll}
  \alpha^{A,B,C}_{x} &:& (A_x\otimes_x B_x)\otimes_x C_x\to A_x\otimes_x(B_x\otimes_x C_x), \\
  \lambda^{A}_{x} &:& I_x\otimes_x A_x\to A_x,\\
  \rho^{A}_{x} &:& A_x\otimes_x I_x\to A_x.
\end{array}
\]
To show that these maps are morphisms of nominal families, by
Remark~\ref{rk:NomFamPhiMap}, we must show that the following diagrams
commute, where $y=\pi\bul x$:
\[
\begin{tikzpicture}
  \node (ltop) at (0,0) {$\pihat((A_x\otimes_x B_x)\otimes_x C_x)$};
  \node (rtop) at (6,0) {$\pihat(A_x\otimes_x(B_x\otimes_x C_x))$};
  \node (lbot) at (0,-1.5) {$(A_y\otimes_y B_y)\otimes_y C_y$}; 
  \node (rbot) at (6,-1.5) {$A_y\otimes_y(B_y\otimes_y C_y)$}; 
  
  \draw [-stealth][black,thick] (ltop) -- (rtop) node[midway,above] {$\pihat(\alpha_{x})$};
  \draw [-stealth][black,thick] (ltop) -- (lbot) node[midway,left] {$\eqfams{(A\otimes_X B)\otimes_X C}{\pi}{x}$};
  \draw [-stealth][black,thick] (lbot) -- (rbot) node[midway,above] {$\alpha_{y}$};
  \draw [-stealth][black,thick] (rtop) -- (rbot) node[midway,right] {$\eqfams{A\otimes_X(B\otimes_X C)}{\pi}{x}$};
\end{tikzpicture}
\]
\[
\begin{tikzpicture}
  \node (ltop) at (0,0) {$\pihat(A_x)$};
  \node (rtop) at (6,0) {$\pihat(I_x\otimes_x A_x)$};
  \node (lbot) at (0,-1.5) {$A_y$}; 
  \node (rbot) at (6,-1.5) {$I_y\otimes_y A_y$}; 
  
  \draw [-stealth][black,thick] (ltop) -- (rtop) node[midway,above] {$\pihat(\lambda_{x})$};
  \draw [-stealth][black,thick] (ltop) -- (lbot) node[midway,left] {$\eqfams{A}{\pi}{x}$};
  \draw [-stealth][black,thick] (lbot) -- (rbot) node[midway,above] {$\lambda_{y}$};
  \draw [-stealth][black,thick] (rtop) -- (rbot) node[midway,right] {$\eqfams{I\otimes_X A}{\pi}{x}$};
\end{tikzpicture}
\quad
\begin{tikzpicture}
  \node (ltop) at (0,0) {$\pihat(A_x)$};
  \node (rtop) at (6,0) {$\pihat(A_x\otimes_x I_x)$};
  \node (lbot) at (0,-1.5) {$A_y$}; 
  \node (rbot) at (6,-1.5) {$A_y\otimes_y I_y$}; 
  
  \draw [-stealth][black,thick] (ltop) -- (rtop) node[midway,above] {$\pihat(\rho_{x})$};
  \draw [-stealth][black,thick] (ltop) -- (lbot) node[midway,left] {$\eqfams{A}{\pi}{x}$};
  \draw [-stealth][black,thick] (lbot) -- (rbot) node[midway,above] {$\rho_{y}$};
  \draw [-stealth][black,thick] (rtop) -- (rbot) node[midway,right] {$\eqfams{A\otimes_X I}{\pi}{x}$};
\end{tikzpicture}
\]
We only prove the diagram for $\alpha$, as the other two are very
similar. Consider the following:
\[
\begin{tikzpicture}
  \node (ltop) at (0,0) {$\pihat((A_x\otimes_x B_x)\otimes_x C_x)$};
  \node (rtop) at (8,0) {$ \pihat(A_x\otimes_x(B_x\otimes_x C_x))$};
  \node (lmid) at (0,-2) {$(\pihat(A_x)\otimes_y \pihat(B_x))\otimes_y \pihat(C_x)$};
  \node (rmid) at (8,-2) {$\pihat(A_x)\otimes_y (\pihat(B_x)\otimes_y \pihat(C_x))$};
  \node (lbot) at (0,-4) {$(A_y\otimes_y B_y)\otimes_y C_y$}; 
  \node (rbot) at (8,-4) {$A_y\otimes_y(B_y\otimes_y C_y)$}; 
  
  \draw [-stealth][black,thick] (ltop) -- (rtop) node[midway,above] {$\pihat(\alpha_{x})$};
  \draw [-stealth][black,thick] (ltop) to[out=200, in=160, looseness=2] node[pos=0.1,left] {$\eqfams{(A\otimes_X B)\otimes_X C}{\pi}{x}$} (lbot);
  \draw [-stealth][black,thick] (ltop) -- (lmid) node[midway,right] {$\mm^3_{\pi}$};
  \draw [-stealth][black,thick] (lmid) -- (lbot) node[midway,right] {\small$(\eqfams{A}{\pi}{x}\otimes_y\eqfams{B}{\pi}{x})\otimes_y\eqfams{C}{\pi}{x}$};
  \draw [-stealth][black,thick] (rtop) to[out=-20, in=20, looseness=2] node[pos=0.1,right=1em] {$\eqfams{A\otimes_X(B\otimes_X C)}{\pi}{x}$}  (rbot);
  \draw [-stealth][black,thick] (rtop) -- (rmid) node[midway,left] {$\mm^{3'}_{\pi}$};
  \draw [-stealth][black,thick] (rmid) -- (rbot) node[midway,left] {\small$\eqfams{A}{\pi}{x}\otimes_y(\eqfams{B}{\pi}{x}\otimes_y\eqfams{C}{\pi}{x})$};
  \draw [-stealth][black,thick] (lmid) -- (rmid) node[midway,above] {$\alpha_{y}$};
  \draw [-stealth][black,thick] (lbot) -- (rbot) node[midway,above] {$\alpha_{y}$};
\end{tikzpicture}
\]
Here, the maps $\mm^3_{\pi}$ and $\mm^{3'}_{\pi}$ are the obvious
morphisms defined by using the $\mm$-maps twice. The parts labelled (1)
commute by definition of $\eqfams{(A\otimes_X B)\otimes_X C}{\pi}{x}$ and
$\eqfams{A\otimes_X (B\otimes_X C)}{\pi}{x}$. The part labelled (2) commutes
because $\mm_{\pi}$ is the structure that makes $\pihat$ into a
monoidal functor. The part labelled (3) commutes by naturality of
$\alpha_y$. 

\section{Proof of Properties (\ref{ax-c}) and (\ref{ax-d}) from Section~\ref{sssec:def-a-magicwand-b}}
\label{app:prop-a-b-2}

\begin{xlem}\label{lem:prop-a-and-b-2}
  The family of maps $\eqfams{C}{\pi}{x}$ defined in
  Section~\ref{sssec:def-a-magicwand-b} satisfies properties
  (\ref{ax-c}) and (\ref{ax-d}) of Remark~\ref{rk:NomFamPhiMap}.
\end{xlem}

\begin{proof}
We must prove properties (\ref{ax-c}) and (\ref{ax-d}) of
Remark~\ref{rk:NomFamPhiMap}. For (\ref{ax-c}), consider any
$\pi_1,\pi_2\in\PI$ and some $x\in X$. We must show that
$\eqfams{C}{\pi_1\pi_2}{x}=\eqfams{C}{\pi_1}{\pi_2\bul x}\circ
\widehat{\pi_1}(\eqfams{C}{\pi_2}{x})$.  Consider the following
diagram. As usual, we have written $y=\pi_2\bul x$ and
$z=\pi_1\pi_2\bul x$.
\[
\scalebox{0.88}{$
\begin{tikzpicture}[scale=0.9]
    \useasboundingbox (-7,-7) rectangle (7,6.9);
    \node (1) at (145:7) {$\widehat{\pi_1\pi_2}(A_x \multimap_x B_x)$};
    \node (2) at (155:7) {$\widehat{\pi_1}(\widehat{\pi_2}(A_x \multimap_x B_x))$};
    \node (3) at (210:4.5) {$\widehat{\pi_1}(\widehat{\pi_2}(A_x)\multimap_{y}\widehat{\pi_2}(B_x))$}; 
    \node (4) at (270:7) {$\widehat{\pi_1}(A_{y}\multimap_{y} B_{y})$};
    
    \node (5) at (30:7) {$A_{z} \multimap_{z} B_{z}$}; 
    \node (6) at (330:4.5) {$\widehat{\pi_1}(A_{y}) \multimap_{z} \widehat{\pi_1}(B_{y})$}; 
    
    \node (7) at (90:4.5) {$\widehat{\pi_1\pi_2}(A_x) \multimap_{z} \widehat{\pi_1\pi_2}(B_x)$}; 
    \node (8) at (0,0) {$\widehat{\pi_1}(\widehat{\pi_2}(A_x))\multimap_{z} \widehat{\pi_1}(\widehat{\pi_2}(B_x))$};
    
    \node (9) at (150:3) {(2)};
    \node (10) at (90:5.75) {(1)};
    \node (11) at (270:3) {(4)};
    \node (12) at (30:3) {(3)};
    \node (13) at (220:5.75) {(1)};
    \node (14) at (320:5.75) {(1)};
    
    \path[every node/.style={font=\sffamily\small}]  
    (2) [-stealth,black,  thick] edge node [right] {$\widehat{\pi_1}(\nn_{\pi_2})$} (3)
    (3) [-stealth,black,  thick] edge node [rotate=-50,above] {$\widehat{\pi_1}((\eqfams{A}{\pi_2}{x})\inv\multimap_{y} \eqfams{B}{\pi_2}{x})$} (4)    
    (4) [-stealth,black,  thick] edge node [midway,left] {$\nn_{\pi_1}$} (6)
    (6) [-stealth,black,  thick] edge node [rotate=69.5,above] {$(\eqfams{A}{\pi_1}{y})\inv\multimap_{z} \eqfams{B}{\pi_1}{y}$} (5)    
    
    (2)  [-stealth,black,  thick] edge [out=250,in=170, looseness=1.1] node [midway,left] {$\widehat{\pi_1}(\eqfams{C}{\pi_2}{x})$} (4)
    (4) [-stealth,black,  thick] edge [out=10, in=290, looseness=1.1] node [midway,right] {$\eqfams{C}{\pi_1}{y}$} (5)
    
    (3) [-stealth,black,  thick,  thick] edge node [midway,above=0.5em] {$\nn_{\pi_1}$} (8)
    (8) [-stealth,black,  thick,  thick] edge node [near end,left] {$\widehat{\pi_1}((\eqfams{A}{\pi_2}{x})\inv) \multimap_{z} \widehat{\pi_1}(\eqfams{B}{\pi_2}{x})$} (6)
    
    (7) [-stealth,black,  thick] edge node [below,rotate=-9] {$(\eqfams{A}{\pi_1\pi_2}{x})\inv\multimap_{z} \eqfams{B}{\pi_1\pi_2}{x}$} (5)
    (1) [-stealth,black,  thick] edge node [midway,below] {$\nn_{\pi_1\pi_2}$} (7)    
    (1) [-stealth,black,  thick, out=50, in=130, looseness=1.1] edge node [midway,above] {$\eqfams{C}{\pi_1\pi_2}{x}$} (5);

    \draw (1) [black,  thick] edge [double distance=1.5pt] (2);
    \draw (7) [black,  thick] edge [double distance=1.5pt] (8);
\end{tikzpicture}
$}
\]
The parts labelled (1) commute by definition of $\varphi^C$. Part (2)
commutes by Remark~\ref{rmk:n-map-comp}. Part (3) commutes because
$\varphi^A$ and $\varphi^B$ belong to nominal families and therefore
satisfy property (\ref{ax-c}). Part (4) commutes by naturality of
$\nn_{\pi_1}$. Therefore, the outside commutes, proving the claim.

For (\ref{ax-d}), assume $\pi\# x$, so that $\pi\bul x = x$. We
must show that $\eqfams{C}{\pi}{x} = \pi \bul (-)$. Consider the
following diagram.
\[
\begin{tikzpicture}
    \node (1) at (-6,0) {$\pihat(A_x\multimap_x B_x)$};         
    \node (2) at (0,0) {$\pihat(A_{x})\multimap_{x} \pihat(B_{ x})$};
    \node (3) at (6,0) {$A_{x}\multimap_{x} B_{x}$};
    \node at (-1,1) {(1)};
    \node at (3,0) {(2)};
    \node at (-1,-1) {(3)};
    
    \path[every node/.style={font=\sffamily\small}]
    (1) [-stealth,black,thick] edge [] node [midway,above] {$\nn_{\pi}$} (2)
    (1) [-stealth,black,  thick] edge [out=-60, in=-120, looseness=0.5] node [midway,below]{$\eqfams{C}{\pi}{x}$} (3) 
    (2) [-stealth,black,  thick] edge [bend right=15] node [midway,below]{$(\eqfams{A}{\pi}{x})\inv\multimap_{x}\eqfams{B}{\pi}{x}$} (3)
    (1) [-stealth,black,  thick] edge [out=60, in=120, looseness=0.5] node [midway,above]{$\DB{\pi}$} (3)
    (2) [-stealth,black,  thick] edge [bend left=15] node [midway,above]{$(\DB{\pi})\inv\multimap_{x} \DB{\pi}$} (3)  ;
    
\end{tikzpicture}
\]
Part (3) commutes by definition of $\varphi^C$. Part (2) commutes
because $\varphi^A$ and $\varphi^B$ each satisfy (\ref{ax-d}). To show
part (1), it is (by the universal property of $\nn$) sufficient to show
that part (1) of the diagram below commutes.

\[
\begin{tikzpicture}
    \node (1) at (-7,4) {$\pihat(A_x\multimap_x B_x) \otimes_x \pihat(A_x)$};
    \node (2) at (7,4) {$\pihat((A_x\multimap_x B_x) \otimes_x A_x)$}; 
    \node (3) at (0,2) {$(A_{x}\multimap_{x} B_{x}) \otimes_x A_{x}$};
    \node (4) at (-4,0) {$(A_{x}\multimap_{x} B_{x}) \otimes_x \pihat(A_x)$};
    \node (5) at (4,0) {$B_{ x}$}; 
    \node (6) at (0,-2) {$(\pihat(A_x) \multimap_{x} B_{x}) \otimes_x \pihat(A_x)$};
    \node (7) at (-7,-4) {$(\pihat(A_{x})\multimap_{x} \pihat(B_{ x})) \otimes_x \pihat(A_x)$};
    \node (8) at (7,-4) {$\pihat(B_{ x})$}; 
    
    \node at (-5,-2) {(1)};
    \node at (-3,2) {(2)};
    \node at (1,3.2) {(3)};
    \node at (0,0) {(4)};
    \node at (3,-3) {(5)};
    \node at (5,2) {(6)};
    
    \path[every node/.style={font=\sffamily\small}]
    (1) [-stealth,black,  thick] edge node [midway,above] {$\mm_{\pi}\inv$} (2)
    (1) [-stealth,black,  thick] edge node [midway,right=0.5em] {$\DB{\pi} \otimes_x \DB{\pi}$} (3)
    (1) [-stealth,black,  thick] edge node [near end,left] {$\DB{\pi} \otimes_x \id$} (4)
    (1) [-stealth,black,  thick] edge node [midway,left] {$\nn_{\pi}\otimes_x \id$} (7)
    (2) [-stealth,black,  thick] edge node [midway,right] {$\pihat(\eps)$} (8)
    (2) [-stealth,black,  thick] edge node [midway,below=0.5em] {$\DB{\pi}$} (3)
    (3) [-stealth,black,  thick] edge node [midway,above] {$\eps$} (5)
    (4) [-stealth,black,  thick] edge node [midway,left=0.5em] {$\id \otimes_x \DB{\pi}$} (3)
    (4) [-stealth,black,  thick] edge node [midway,left=0.5em] {$(\DB{\pi}\multimap_{x} \id) \otimes_x \id$} (6)
    (8) [-stealth,black,  thick] edge node [near end,right] {$\DB{\pi}$} (5)
    (7) [-stealth,black,  thick] edge node [midway,right=0.5em] {$(\id \multimap_{x} \DB{\pi}) \otimes_x \id$} (6)
    (6) [-stealth,black,  thick] edge node [midway,below] {$\eps$}  (5)
    (7) [-stealth,black,  thick] edge node [midway,below] {$\eps$}  (8);
\end{tikzpicture}
\]
The outer square commutes by definition of the $\nn$ map, part (2) commutes by composition, part (3) commutes Lemma~\ref{lem:m-pi}, part (4) commutes by dinaturality of $\eps$, part (5) commutes by naturality of $\eps$, and part (6) commutes by naturality of $(\DB{\pi})\inv$. Therefore, part (1) commutes, as well.
\end{proof}

\end{document}